\documentclass[12pt]{amsart}
\usepackage{amssymb}

\usepackage{enumerate}

\usepackage{diagrams}
\diagramstyle[scriptlabels,height=8mm,width=8mm]

\def\bdi{\begin{diagram}}
\def\edi{\end{diagram}}

\theoremstyle{plain}

\newtheorem{thm}{Theorem}[section]
\newtheorem{cor}[thm]{Corollary}
\newtheorem{lem}[thm]{Lemma}
\newtheorem{prop}[thm]{Proposition}

\theoremstyle{definition}
\newtheorem{defi}[thm]{Definition}
\newtheorem{defis}[thm]{Definitions}
\newtheorem{conj}[thm]{Problem}
\newtheorem{conv}[thm]{Convention}
\newtheorem{nota}[thm]{Notation}
\newtheorem{rem}[thm]{Remark}
\newtheorem{rems}[thm]{Remarks}
\newtheorem{exa}[thm]{Example}
\newtheorem{exas}[thm]{Examples}
\newtheorem{sit}[thm]{}

\newcommand{\brem}{\begin{rem}}
\newcommand{\brems}{\begin{rems}}
\newcommand{\erem}{\end{rem}}
\newcommand{\erems}{\end{rems}}
\newcommand{\bexa}{\begin{exa}}
\newcommand{\bexas}{\begin{exas}}
\newcommand{\eexa}{\end{exa}}
\newcommand{\eexas}{\end{exas}}
\newcommand{\bdefi}{\begin{defi}}
\newcommand{\edefi}{\end{defi}}
\newcommand{\bdefis}{\begin{defis}}
\newcommand{\edefis}{\end{defis}}
\newcommand{\bcor}{\begin{cor}}
\newcommand{\ecor}{\end{cor}}
\newcommand{\blem}{\begin{lem}}
\newcommand{\elem}{\end{lem}}
\newcommand{\bconv}{\begin{conv}}
\newcommand{\econv}{\end{conv}}
\newcommand{\bconj}{\begin{conj}}
\newcommand{\econj}{\end{conj}}
\newcommand{\bprop}{\begin{prop}}
\newcommand{\eprop}{\end{prop}}
\newcommand{\bthm}{\begin{thm}}
\newcommand{\ethm}{\end{thm}}
\newcommand{\bnota}{\begin{nota}}
\newcommand{\enota}{\end{nota}}
\newcommand{\bsit}{\begin{sit}}
\newcommand{\esit}{\end{sit}}
\newcommand{\be}{\begin{equation}}
\newcommand{\ee}{\end{equation}}
\newcommand{\bproof}{\begin{proof}}
\newcommand{\eproof}{\end{proof}}

\def\ba{\begin{array}}
\def\ea{\end{array}}
\def\bea{\begin{eqnarray}}
\def\eea{\end{eqnarray}}

\def\bnum{\begin{enumerate}}
\def\enum{\end{enumerate}}
\newcommand{\no}{\noindent}

\def\lto{\longrightarrow}
\def\hto{\hookrightarrow}

\def\ext{{\rm ext}}
\def\sto{\rightsquigarrow}
\def\and{\quad \mbox{and}\quad}

\newcommand{\Sing}{\operatorname{Sing}}
\newcommand{\Spec}{\operatorname{Spec}}

\newcommand{\Pic}{\operatorname{Pic}}

\newcommand{\id}{\operatorname{id}}

\newcommand{\Aut}{{\operatorname{Aut}}}
\newcommand{\ML}{{\operatorname{ML}}}

\newcommand{\reg}{{\operatorname{reg}}}

\newcommand{\Ext}{{\operatorname{Ext}}}
\newcommand{\Hom}{{\operatorname{Hom}}}

\newcommand{\cHom}{{\mathcal Hom}}

\newcommand{\GL}{{\bf {GL}}}

\newcommand{\Reg}{{\operatorname{Reg}}}

\def\loc{{\rm{loc}}}

\def\fm{{\mathfrak m}}

\def\fC{{\mathfrak C}}

\def\fF{{\mathfrak F}}

\def\fM{{\mathfrak M}}

\def\fR{{\mathfrak R}}

\def\cA{{\mathcal A}}
\def\cB{{\mathcal B}}
\def\cC{{\mathcal C}}
\def\cD{{\mathcal D}}
\def\cE{{\mathcal E}}
\def\cF{{\mathcal F}}
\def\cG{{\mathcal G}}

\def\cL{{\mathcal L}}
\def\cM{{\mathcal M}}

\def\cO{{\mathcal O}}
\def\cP{{\mathcal P}}

\def\cS{{\mathcal S}}

\def\cV{{\mathcal V}}

\def\cX{{\mathcal X}}

\def\cZ{{\mathcal Z}}

\renewcommand{\AA}{{\mathbb A}}

\newcommand{\GG}{{\mathbb G}}
\newcommand{\HH}{{\mathbb H}}
\newcommand{\PP}{{\mathbb P}}
\newcommand{\RR}{{\mathbb R}}
\newcommand{\CC}{{\mathbb C}}
\newcommand{\QQ}{{\mathbb Q}}
\newcommand{\ZZ}{{\mathbb Z}}

\newcommand{\kk}{{\Bbbk}}

\def\tcV{{\tilde{\mathcal V}}}
\def\bcV{{\bar{\mathcal V}}}

\def\bB{{\bar B}}
\def\bC{{\bar C}}
\def\bD{{\bar D}}

\def\bV{{\bar V}}
\def\bX{{\bar X}}

\def\tB{{\tilde B}}
\def\tV{{\tilde V}}

\newcommand{\la}{\label}

\renewcommand{\phi}{\varphi}
\renewcommand{\rho}{\varrho}

\newcommand{\G}{{\Gamma}}

\def\ff{{\bf F}}
\def\An{{\bf ASch_\kk}}
\def\sets{{\bf sets}}

\def\reg{{reg}}

\newcommand{\nlin}{\unitlength1mm\begin{picture}(0,9.25)
                     \put(0,0.75){\line(0,1){8.5}}
                    \end{picture}}

\newcommand{\vlin}[1]{\hspace{0.75mm}\unitlength1mm\begin{picture}
(#1,0)
                     \put(0,0){\line(1,0){#1}}
                    \end{picture}\hspace{0.75mm}\rule[-3mm]{0mm}{4mm}}

\def\llin{\vlin{11.5}}
\newcommand{\lin}{\vlin{8.5}}

\newcommand{\co}[1]{\unitlength1mm\begin{picture}(0,8)
  \put(0,0){\circle{1.5}}
  \put(0,3){\makebox(0,5)[b]{$#1$}}
                    \end{picture}}
\newcommand{\mybox}{\unitlength1mm\begin{picture}(0,1.5)
  \put(-0.75,-0.75){\line(0,1){1.5}}
  \put(-0.75,-0.75){\line(1,0){1.5}}
  \put(0.75,0.75){\line(0,-1){1.5}}
  \put(0.75,0.75){\line(-1,0){1.5}}
  \end{picture}}

\newcommand{\xbox}{\unitlength1mm\begin{picture}(0,1.5)
  \put(0,0){$\mybox$}
  \put(-0.75,0){\line(1,0){1.5}}
  \put(0,-0.75){\line(0,1){1.5}}
  \end{picture}}

\newcommand{\cou}[2]{\unitlength1mm\begin{picture}(0,8)
  \put(0,0){\circle{1.5}}
  \put(0,3){\makebox(0,5)[b]{$#1$}}
  \put(0,-7){\makebox(0,4)[t]{$#2$}}
    \end{picture}
    \rule[-7mm]{0mm}{7mm}}

\newcommand{\crl}[2]{\unitlength1mm\begin{picture}(0,8)
  \put(0,0){\circle{1.5}}
  \put(-5,0){\makebox(0,5)[b]{$#1$}}
 \put(5,0){\makebox(0,5)[b]{$#2$}}
    \end{picture}
    \rule[-7mm]{0mm}{7mm}}

\newcommand{\cshiftup}[2]{\unitlength1mm\begin{picture}(0,9.25)
                     \put(0,10){\crl{#1}{#2}}
                    \end{picture}}

\newcommand{\xbrl}[2]{\unitlength1mm\begin{picture}(0,8)
  \put(0,0){\xbox}
  \put(-5,0){\makebox(0,5)[b]{$#1$}}
 \put(5,0){\makebox(0,5)[b]{$#2$}}
    \end{picture}
    \rule[-7mm]{0mm}{7mm}}

\newcommand{\xbshiftup}[2]{\unitlength1mm\begin{picture}(0,9.25)
                     \put(0,10){\xbrl{#1}{#2}}
                    \end{picture}}

\title{Deformation equivalence of affine ruled surfaces}

\author{Hubert Flenner}
\address{Fakult\"at f\"ur Mathematik,
Ruhr Universit\"at Bochum, Geb.\ NA 2/72, Universit\"ats\-str.\
150, 44780 Bochum, Germany}
\email{Hubert.Flenner@ruhr-uni-bochum.de}

\author{Shulim Kaliman}
\address{Department of Mathematics,
University of Miami, Coral Gables, FL  33124, U.S.A.}
\email{kaliman@math.miami.edu}

\author{Mikhail Zaidenberg}
\address{Universit\'e
Grenoble I, Institut Fourier, UMR 5582 CNRS-UJF, BP 74, 38402
St.\ Martin d'H\`eres c\'edex, France}
\email{zaidenbe@ujf-grenoble.fr}

\thanks{
{\bf Acknowledgements:} This research was done during a visit of
the first and the second authors at the  Fourier Institute,
Grenoble, and of all three authors at
the Max-Planck-Institute of Mathematics, Bonn.
They thank these institutions for the
generous support and excellent working conditions.}

\thanks{
\mbox{\hspace{11pt}}{\it 2010 Mathematics Subject
Classification}:
14R05, 14R25, 14J10, 14D22.\\
\mbox{\hspace{11pt}}{\it Key words}: affine ruled surface,
$\AA^1$-fibration, moduli space, deformation equivalence}

\date{\today}

\begin{document}

\begin{abstract} 
A smooth family $\phi:\cV\to S$ of surfaces will be called {\em completable} if there is a logarithmic deformation $(\bcV,\cD)$ over $S$ so that $\cV=\bcV\backslash \cD$. 
Two smooth surfaces $V$ and $V'$ are said to be deformations of each other if there is a completable flat family $\cV\to S$ of smooth surfaces over a connected base so that $V$ and $V'$ are fibers over suitable points $s,s'\in S$. This relation generates an equivalence relation called {\em deformation equivalence}. In this paper we give a complete combinatorial description of this relation in the case of affine ruled surfaces, which by definition are surfaces that admit an affine ruling $V\to B$ over an affine base with possibly degenerate fibers. In particular we construct complete families of such affine ruled surfaces. In a few particular cases we can  also deduce the existence of
a coarse moduli space.
\end{abstract}

\maketitle

\tableofcontents

\section*{Introduction}\label{intro}

In classifying algebraic objects like varieties or vector bundles  of a certain class one usually tries to find discrete invariants such that the objects sharing a fixed set of these invariants, form a moduli space. The model case is here the moduli space of smooth complete curves,  which is known to be a $3g-3$-dimensional irreducible variety, see \cite{DeMu}. 
Another classical case is the Hilbert scheme $\HH_P$ parameterizing subschemes of the projective space $\PP^n$ with fixed Hilbert polynomial $P$. By a classical result of Hartshorne \cite{Ha1} $\HH_P$ is connected. On the other hand it is already an unsolved problem to determine the connected components of the Hilbert scheme of locally Cohen-Macaulay curves of degree $d$ and genus $g$, see \cite{Ha3} and the references therein for partial results.

In general it may happen that for a class of varieties one cannot expect a reasonable moduli space. A typical obstruction for its existence are usually $\PP^1$-fibrations or, in the case of open varieties, $\AA^1$-fibrations. Instead of studying in such cases the connected components of the moduli space it is convenient to consider the relation of deformation equivalence as introduced in  \cite{Ca}.
Given a deformation $f:\cX\to S$ with all fibers in a class of varieties $\cC$ we consider the relations $f^{-1}(s)\sim f^{-1}(t)$ if and only if $s,t\in S$ belong to the same connected component of $S$. This generates an equivalence relation on $\cC$ called {\em deformation equivalence}. 
 
In this paper we study this equivalence relation for the class of  normal  affine surfaces $V$ that can be equipped with an {\em affine ruling}. By this we mean a morphism $\pi:V\to B$ onto a smooth {\em affine} curve with the general fiber isomorphic to $\AA^1$. In this case there is no moduli space available except for a few special cases; see sect.\ 5.
Our main result is that nevertheless  one can characterize  deformation equivalence completely.  The deformations which we consider here are deformations of the open surface $X$ which can be extended to a logarithmic deformation of a completion of $X$ in the sense of \cite{Ka}.

The main tool for this characterization is the so called  {\em normalized extended graph}. For smooth surfaces it can be described as follows. Given an affine ruling $\pi:V\to B$  we can extend $\pi$ to a $\PP^1$-fibration $\bar\pi:\bX\to\bB$ over the normal completion $\bB$ of $B$. Performing suitable blowups and blowdowns along the boundary $D:=\bX\backslash X$   of $X$ one can transform $D$ or, equivalently,  its weighted dual graph $\Gamma_D$ into standard form, see \cite{DaGi} or \cite{FKZ1}. The {\em extended divisor} $D_\ext$ is then the union of $D$ with all  singular fibers of $\bar\pi$. 
This extended divisor is  not  in general an invariant of the surface alone, but depends on the choice of the affine ruling.    
The components of $D_\ext-D$ are called the {\em feathers}. For such feathers one has the notion of a {\em mother component} (see Definition \ref{mother} for details). A feather is not necessarily linked to its mother component.  Attaching the feathers to their mother components we obtain from  the dual graph $\G_\ext$ of $D_\ext$ the {\em normalized extended graph} $N(\Gamma_\ext)$. With this terminology our main result is as follows.

\bthm
Two affine ruled surfaces are deformation equivalent if and only if they share the same normalized extended graph unless $D$ is a zigzag, i.e.\ a  linear chain of rational curves. In the latter case two surfaces are deformation equivalent if and only if the normalized extended graphs are equal, possibly after reversing one of them (see sect.\ 3 for details). 
\ethm

This theorem enables us to construct in section 5 a coarse moduli space for the class of special Gizatullin surface that was studied in \cite{FKZ4}.

Let us mention a few special cases.  First
of all, the Ramanujam theorem \cite{Ra} states that every contractible,
smooth surface  over $\CC$ with a trivial fundamental group at
infinity (in particular, every smooth
surface over $\CC$ homeomorphic to $\RR^4$)
is isomorphic to the affine plane $\AA^2=\AA_{\CC}^2$. Thus the
deformation equivalence class for such surfaces is reduced to a single point. The normalized extended graph consists in this case just of two 0-vertices joined by an edge. 

A generalization of this theorem due to Gurjar and Shastri
\cite{GuSh} says that every normal contractible surface with a
finite fundamental group at infinity is isomorphic to $\AA^2/G$,
where $G$ is a finite subgroup of $\GL_2(\CC)$ acting freely on $\AA^2\backslash \{ 0 \}$.
Fixing $G$, the deformation equivalence class for such surfaces again consists of a single point.

Yet another result of this type is provided by
the Danilov-Gizatullin Isomorphism Theorem
\cite{DaGi} (see also \cite{CNR, FKZ5}).
Recall that a {\em Danilov-Gizatullin
surface} is the complement to an ample divisor $H$ in some Hirzebruch surface.  Such a  surface can be completed by the zigzag $[[0,0,(2)_n]]$, where $n+1=H^2$. The
Danilov-Gizatullin theorem proves that the deformation equivalence class for such a surface again consists of a single point, once we fix the length of the boundary
zigzag.

The paper is organized as follows. In section \ref{prelim} we
recall the notion of standard completion and standard dual
graph of an affine ruled surface, and in subsection  \ref{extdivi}
that of an extended divisor. The rest of  section  \ref{famafsur}
is devoted to developing such notions in the relative case. 
More specifically, in Theorem \ref{1.6} we show the existence of standard completions in completable families,  and in the   factorization Lemma \ref{factlem} we shed some light onto the structure of extended divisors in the relative case. 
Central  are here the notions of {\em completable families} and {\em resolvable families}, see Definitions \ref{NCcomplet} and \ref{def-res-comp}. The
normalized extended graph is introduced in section \ref{defequiv}.
Our main result Theorem \ref{main1} characterizes the
deformation equivalence of affine ruled surfaces in terms of these
graphs. The proof of this theorem starts in section
\ref{defequiv} and is completed in the next section 4, where we
construct a versal deformation space of affine ruled surfaces. In
section \ref{modspecgiz} we apply our previous results in order to
construct the moduli space of special Gizatullin surfaces, a
subclass of the class of all Gizatullin surfaces studied in detail
in our previous paper \cite{FKZ4}.

All varieties in this paper are assumed to be defined over an algebraically closed field $\kk$. By a {\em surface} we mean a  connected algebraic scheme over $\kk$ such that all its irreducible
components are of dimension 2.

\section{Preliminaries on standard completions}\la{prelim}

Let us recall the notions of a (semi-)standard zigzag and
(semi-)standard graph. 

\bsit\label{zigz} Let $X$ be a complete normal algebraic surface,
and let $D$ be an SNC (i.e.\ a simple normal crossing) divisor $D$
contained in the smooth part $X_{\rm reg}$ of
$X$. We say that $D$ is a {\it zigzag} if all irreducible components of $D$ are
rational and the dual graph $\G_D$ of $D$ is linear\footnote{By abuse
of notation,
we often denote an SNC divisor and its dual graph by the same
letter.}.
We abbreviate a chain of curves $C_0, C_{1},\ldots, C_n$ of weights
$w_0,\ldots, w_n$ by $[[w_0,\ldots, w_n]]$. We also write
$[[\ldots, (w)_k,\ldots ]]$ if a weight $w$ occurs at $k$
consecutive places. \esit

\bsit\label{stzi} A zigzag $D$ is called {\em standard}
if it is one of the chains
\be \label{standardzigzag}
[[(0)_i]]\,,i\le 3,\,\,\mbox{or}\,\,
  [[(0)_{2i},w_2,\ldots,w_n]],\,\, \mbox{where}\,\, i\in \{0,1\},
\,\, n\ge 2 \,\, \mbox{and}\,\, w_j\le-2\,\,\,\forall j.
\ee
A linear chain $\Gamma$ is said to be {\em semi-standard} 
or {\em $w_1$-standard} if it is either standard or one of
\be \label{sstandard}
 [[0,w_1,w_2,\ldots,w_n]],\quad [[0,w_1,0]]\,, \,\, \mbox{where}
\,\, n\ge 1,\,\,w_1\in \ZZ,\,\,\mbox{and}\,\, w_j\le-2\,\,\,\forall
j\ge 2\,. 
\ee
A {\em circular
graph} is a connected graph with all vertices of degree 2. Such a
weighted graph will be denoted by $((w_0,\ldots,w_n))$. A circular
graph is called {\em standard} if it is one of
\begin{equation}\label{cstandard}
((w_1,\ldots, w_n)), \quad\quad ((0,0, w_1,\ldots, w_n)),
\quad\quad ((0_l, w)), \quad\mbox{or}\quad ((0,0, -1,-1)) ,
\end{equation}
where $w_1,\ldots,w_n\le -2$, $n>0$, $0\le l\le3$ and  $w\le
0$.\footnote{We list only those graphs for which the intersection
matrix has at most one positive eigenvalue.
Actually we do not use circular graphs in this paper, however without
 these graphs
the general notion of a standard graph in \ref{stgraph} below would be incomplete.}

By an {\em
inner elementary transformation} of a weighted graph we mean
blowing up at an edge incident to a $0$-vertex of degree $2$ and
blowing down the image of this vertex. By a sequence of inner
elementary transformations we can successively move the pair of
zeros in the standard zigzag $[[0,0,w_2,\ldots,w_n]]$ to the
right:
$$
[[0,0,w_2,\ldots,w_n]]\sto [[w_2, 0,0,w_3,\ldots,w_n]] \sto \ldots
\sto [[w_2,\ldots,w_n, 0,0]]\,.
$$
This yields the {\em reversion} 
\be\la{reverse}
D=[[0,0,w_2,\ldots,w_n]]\rightsquigarrow
[[0,0,w_n,\ldots,w_2]]=:D^\vee\, 
\ee 
(see 1.4 in \cite{FKZ3}).

An {\em outer elementary transformation} 
consists in blowing up at
a $0$-vertex of degree $\le 1$ and 
blowing down the image of this vertex. 
A birational inner elementary transformation on a surface
is rigid i.e., uniquely determined by the 
associated combinatorial transformation of the dual graph, 
whereas an outer one depends on
the choice of the center of blowup.
\esit

\bsit\la{NC-completion}
In the sequel there is also the need to consider {\em NC divisors} $D$ on a surface $X$. By this we mean  that $D\subseteq X_\reg$ and that  the singularities of $D$ are ordinary double
points; these are given in the local ring $\cO_{X,p}$ by an equation $xy=0$, where $x,y\in \fm\backslash\fm^2$.
In particular,  two different components
meet in smooth points transversally, while the intersection of
three different components is empty. Thus $D$ is an SNC divisor  if and only if it is an NC divisor and all its irreducible
components are smooth. The dual graph $\G_D$ of an NC divisor $D$ has loops which correspond to the singular points of the components.
Vice versa, if the dual graph $\G_D$ of an NC divisor $D$ has
no loop then $D$ is SNC. In particular, $D$ is SNC if $\G_D$ is a
tree. 

An {\em NC completion} $(\bV,D)$ of a surface
$V$ consists of a complete surface $\bV$ and an NC divisor $D$
on $\bV$ such that $V=\bV\setminus D$. 
\esit

For instance,
a plane nodal cubic
is an NC divisor, which is not SNC. Its dual graph
consists of one vertex of weight $3$ and a loop.

\bdefi\la{stgraph} Let $\Gamma$ be the dual graph of an NC divisor
$D$. We use the following notations. 
\bnum[(1)]
\item  $B=B(\G)$ is the set of branching points  of $\G$.
\item $S=S(\G)$ is
the set of vertices corresponding to non-rational components of $D$. 
\item Following \cite{FKZ1} a connected component of
$\G- (B\cup S)$ is called a {\em segment} of $\Gamma$. 
\item A
segment will be called  {\em outer} if 
it contains an {\em extremal}  (or {\em end})
vertex of $\G$ i.e., a vertex  of degree 1.
\enum
Thus an outer segment is either the whole graph $\G$ and $\G$  is linear, or it is connected to exactly one vertex of
$B\cup S$. The dual graph $\G$ of
an NC divisor $D$ on an algebraic surface $\bar V$ (and also $D$ itself)
will be called {\it (semi-)standard} if the following hold. 

\bnum[(i)]
\item  All segments of $\Gamma$ are (semi-)standard;

\item 
If a segment is outer and contains a vertex $v$ of weight 
$0$ then it also has an extremal vertex of weight $0$.
For a standard graph we require additionally 
that the neighbor in $\Gamma$ of
every extremal vertex of weight $0$ is as well a zero vertex.
\enum
An NC completion $(\bar V,D)$ of an open  surface $V$ is
called {\em (semi-)standard} if so is $D$.
\edefi

These notions differ from those in \cite[Definition 2.13]
{FKZ1}, where condition (ii) is absent.

\brem\la{1.4}
(1) Every normal surface $V$ has a standard 
NC completion $(\bV,D)$. For applying 
\cite[Theorem 2.15(b)]{FKZ1} every normal surface has an NC completion $(\bV,D)$ such that $\G_D$
satisfies (i) in \ref{stgraph}. By further elementary
transformations we can achieve that also (ii) holds. 

(2) The dual  graph $\G_D$ of the boundary divisor of a standard NC-completion $(V,D)$ is unique up to elementary transformations as follows from \cite[Theorem 3.1]{FKZ1}. 

(3) A Gizatullin surface is a normal affine surface $V$ with a completion $(\bV,D)$, where $D$ is a standard zigzag. Reversing $D$ by a sequence of inner elementary transformations performed on $(\bV, D)$ we obtain a
new completion $(\bV^\vee, D^\vee)$, which is called the {\em
reversed standard completion}. It is uniquely determined by $(\bV,D)$. 
\erem

In the presence of an affine ruling we have more precise informations.  

\blem\label{lemstd}
Let $V$ be a normal affine surface.
Given an  affine ruling $\pi:V\to B$ over a smooth affine curve 
$B$ the following hold.
\bnum[(a)]
\item There is a standard SNC completion $(\bV,D)$ such that $\pi$ extends to a $\PP^1$-fibration $\bar\pi:\bX\to \bB$ over the smooth completion $\bB$ of $B$. There is a unique curve, say, $C_1$ in $D$, which is a section of $\bar\pi$.

\item With $(\bV, D)$ as in (a), $\G_D$ is a tree, and $D-C_1$ has only rational components.

\item A completion as in (a) is unique up to elementary transformations in extremal zero vertices of $\G_D$. These extremal zero vertices correspond to components of $D-C_1$ that are full fibers of $\bar\pi$.
\enum
\elem 

\bproof
(a) and (b): Let $(\bV,D)$ be an SNC completion of $V$. Blowing up $(\bV,D)$ suitably we may assume that $\pi$ extends to a regular map $\bar\pi:\bV\to \bB$. There is a unique horizontal component, say, $C_1$ of $D$ which is a
section of $\bar\pi$. 
The surface $V$ being affine $D$ is
connected. Every fiber of $\bar\pi$ is a tree of rational curves 
and so $D$ is as well a tree with rational components except possibly for $C_1$. 

Any branch $\cB$ of $D$ at $C_1$ (i.e.\ connected
component of $D- C_1$) is contained in a fiber of $\bar\pi$.
If the intersection form of $\cB$ is not negatively definite, then by Zariski's lemma \cite[4.3]{FKZ1} $\cB$ coincides with the entire fiber of
$\bar\pi$ and so can be contracted to a single $0$-curve. Thus after contracting at most linear $(-1)$-curves of $D$ contained in the fibers we
may suppose that every branch of $\G_{D}$ at $C_1$
either is minimal and negative definite or is a single extremal
zero vertex. Since $B$ is affine, extremal zero vertices do exist. Performing
outer elementary transformations at such a vertex
we can achieve that $(C_1)^2=0$. Thus conditions (i) and (ii) of
Definition \ref{stgraph} are fulfilled and so $\G_{D}$ is
standard.

(c) Let $(\bV', D')$ be a second SNC completion satisfying (a). We can find a common domination of $(\bV, D)$ and $(\bV', D')$ by a SNC completion $(\bV'',D'')$ of $V$. We may suppose that there is no $(-1)$-curve in $\bV''$ that is contracted in both $\bV$ and $\bV'$. If $C\ne C_1$ is a curve in $D$, which is not a full fiber, then its proper transform in $\bV''$ has self-intersection $\le -2$.
Hence if a $(-1)$-curve $E$ in $\bV''$ is contracted in $\bV'$, then it is necessarily the proper transform of a component of $D$ which is a full fibers of $\bar\pi$. Thus the indeterminacy locus of the map $\bdi \bV&\rDotsto& \bV'\edi$ is contained in the union of full fibers, say, $F_1\cup\ldots\cup F_a$ of $\bar\pi$  that are contained in $D$. 
Now it is a standard fact that $(\bV',D)$ is obtained  via a sequence of elementary transformations along $F_1\cup\ldots\cup F_a$.
\eproof

\brems\label{cplusa} 
(1) Given a normal affine surface $V$ with a semi-standard $NC$
completion $(\bV, D)$,  the following hold. 
\bnum[(a)]
\item If $\Gamma_D$ contains an extremal 0-vertex, then $V$ is affine ruled. Indeed, every extremal 0-vertex $C_0$ of $\G_D$ induces a $\PP^1$-fibration $\bar\pi:\bV\to \bB$ onto a smooth complete curve $\bB$ so that $C_0$ is one of the fibers, see
e.g.\ \cite[Chapt.\ V, Proposition 4.3]{BHPV}. In particular $\bar\pi$ restricts to an affine ruling on $V$.

\item 
Conversely, if $V$ carries an affine ruling then $\G_D$ has extremal 0-vertices. 
For by Remark \ref{1.4}(1)  the given semi-standard 
completion is obtained by elementary transformations 
from a standard one with extremal 0-vertices as in 
Lemma \ref{lemstd}. Applying \cite[Corollary 3.33']{FKZ1}, 
(b) follows. 
\enum

(2) Recall that a normal affine surface different from
$\AA^1\times\AA^1_*$ \footnote{As before here
$\AA^1_*=\AA^1\setminus\{0\}$.} is Gizatullin if and only if it
admits two distinct affine rulings \cite{Gi, Du1}. 
Thus the affine ruling of a non-Gizatullin surface is unique up to an isomorphism of the base.
Actually the base of the canonical affine ruling $V\to
B$ is equal to  $B:=\Spec \ML(V)$, where $ \ML(V)\subseteq \cO(V)$ is the Makar-Limanov invariant of $V$, i.e.\  the common kernel of  nonzero locally nilpotent derivations on $\cO(V)$. Clearly, the ruling $V\to B$ is induced by the embedding $\ML(V)\hookrightarrow \cO(V)$.

(3) We note that by Remark (1) any standard completion $(\bV,D)$ of an affine ruled surface $V$ yields an affine ruling $\pi:V\to B$ so that $(\bV,D)$ is a standard completion associated to $\pi$ as in Lemma \ref{lemstd}. 

(4) In the situation of Lemma \ref{lemstd}, performing
elementary transformations in extremal 0-vertices one can replace a standard completion by a semi-standard one where $(C_1)^2$ is any given number. This will be useful in the sequel.

\erems

\blem\la{bd-div}
Given an affine ruled surface $V$ with standard completion 
$(\bV,D)$ the following hold.
\bnum[(a)] 
\item If $V$ is Gizatullin then up to 
reversion (\ref{reverse}) the  dual graph $\Gamma_D$ 
of $D$ does not depend on the choice of such a completion.

\item  If $V$ is non-Gizatullin then the divisor $D$ and its dual graph $\G_D$ are unique up to elementary transformations at extremal 0-vertices of $\G_D$. 
\enum
\elem

\bproof
(a) is shown in \cite[Theorem 2.1 ]{DaGi}, see also 
\cite[Corollary 3.32]{FKZ1}.

(b) With $\pi:V\to B$ the affine ruling of $V$ we let $(\bV,D)$ be a standard completion of $V$ as constructed in Lemma \ref{lemstd}(a). By Remark \ref{cplusa}(1) the dual graph $\G_{D'}$ of any other standard completion $(\bV',D')$ has again an extremal 0-vertex that induces a $\PP^1$-fibration  $\bar\pi':\bV'\to \bB'$. Since $V$ is non-Gizatullin, up to a suitable isomorphism $\bB\cong\bB'$ the map $\bar\pi'$ restricts again to $\pi$ on the open part $V$. Thus $(\bV',D')$ is also a standard completion as in Lemma \ref{lemstd}(a). Now the assertion follows form part (c) of that Lemma. 
\eproof

\section{Families of affine surfaces}\label{famafsur}

\subsection{Families of standard completions}\label{famstacom}
By a {\em family} (sometimes called a {\em flat family}) we mean a flat morphism $\phi:\cV\to S$ of algebraic $\kk$-schemes. We call it a family of normal or affine surfaces if every fiber $\cV(s)=\phi^{-1}(s)$, $s\in S$, has this property.

\bdefi\label{reldiv} Let $\phi:\bar \cV\to S$ be a proper  family of surfaces over an algebraic $\kk$-scheme $S$. 
A Cartier divisor
$\cD$ in $\cV$ will be called a {\em relative NC divisor} or a {\em family of NC-divisors} if
$\cD$ is proper over $S$ and for every point $p\in \cD$ either
$\phi|\cD$ is smooth at $p$, or in a suitable neighborhood of $p$
we have $\cD=\cD'\cup\cD''$, where $\cD'$, $\cD''$ and the scheme theoretic intersection $\cD'\cap
\cD''$ are smooth over $S$.

Every fiber $\cD(s)$, $s\in S$, of a relative NC-divisor is an NC divisor. Its
dual graph $\G(s)$ is locally constant in $S$ with respect to the \'etale topology.
We say that a relative NC divisor $\cD$ {\em has constant dual graph} $\G$ if for every
irreducible component $S'$ of $S$ and all $s\in S'$ the
irreducible components of both $\cD'=\cD|S'$ and $\cD(s)$ are in
one-to-one correspondence with the vertices of $\G$.
\edefi

\brem\la{rem2.2}
(1) A relative NC-divisor is contained in the set $\Reg(\cV/S)$ of points in $\cV$ in which $\phi$ is smooth.

(2)  If $\Sing(\cD/S)$ is non-empty then
$\Sing(\cD/S)\to S$ is an unramified covering.
It may
happen that the irreducible components of $\cD$ and $\cD(s)$ are
not in one-to-one correspondence even if the base $S$ is
connected (see Example \ref{jump} below).
In such cases the family $(\G(s))_{s\in S}$ has non-trivial
monodromy, and $\Sing(\cD/S)\to S$ is a non-trivial
covering of $S$.

(3) If $\phi$ is smooth and all irreducible components of $\cD$ are smooth then a relative NC-divisor $\cD$ amounts locally in $S$ to a logarithmic deformation of the fiber $(\cV(s),\cD(s))$ in the sense of \cite{Ka}. 
\erem

The following definition is central in our considerations. 

\bdefi\la{NCcomplet}
A flat family $\phi:\cV\to S$ of normal surfaces will be called
{\em completable} if there is a proper smooth family of surfaces
$\bar\phi: \bar\cV \to S$ together with an open embedding $\cV\hto
\bar\cV$ such that $\cD=\bar \cV\backslash \cV\to S$ is a
relative NC divisor with constant dual graph $\G$. We call such a pair
$(\bar\cV, \cD)$ an {\em NC completion} of $\phi$. This NC completion is called {\em (semi-)standard} if so is $\G$.
\edefi

The main result in this subsection is the
following relative version of the existence of standard NC completions, see Remark \ref{1.4}(a).

\bthm\la{1.6} 
Let $\phi:\cV\to S$ be a family of normal 
surfaces admitting an NC completion $(\bar\cV,\cD)$. Then locally with respect to the \'etale topology there exists a standard NC completion $(\bar\cV ,\cD)$ of
$\phi$.
\ethm

To deduce this result we need some preparations. 
The following  well known lemma  is  an easy consequence of a result due   to Grothendieck; in the analytic setup this is a special case of Kodaira's theorem on the stability of certain submanifolds in deformations \cite{Ko}. 

\blem\la{deflem}  
Let $\phi: \cV\to S$ be a flat family of normal
surfaces over $S$. If the regular part of a special fiber
$V=\cV(s)$ contains a $(-1)$-curve $E$ then over some \'etale neighbourhood of $s$ 
there is a unique family of $(-1)$-curves
$\cE\subseteq\cV$ with $E=\cE(s)$.
\elem

\bproof
The curve $E$ can be considered as a point $[E]$ in the relative Hilbert scheme $\HH_{X/S}$. Applying the smoothness criterion of Grothendieck \cite[Corollaire 5.4]{Gro1} the morphism $\HH_{X/S}\to S$ is \'etale at $[E]$. Hence on a sufficiently small open neighbourhood $U$ of $[E]$ in $\HH_{X/S}$ the map $U\to S$ is \'etale. Now the universal subspace of $X\times_S U$ yields the desired family of $(-1)$-curves.
\eproof

%

In the next result we show that in the relative case
one can perform blowups and
blowdowns just as in the absolute case. In order to keep the formulation short, let us say that a divisor $\cE\subseteq \bar\cV$ is a {\em relative $(-1)$-curve over $S$} if it is smooth over $S$ and restricts to a $(-1)$-curve in every fiber.  

\blem\la{simbd} 
Let $\phi:\bar \cV\to S$ be a proper flat family of
surfaces, and let $\cD$ be a relative NC divisor on $\bar \cV$
over $S$ with constant dual graph $\G$. Then the following hold.
\bnum[(a)]
\item  If $\cD_i\subseteq \cD$ is a relative $(-1)$-curve, which
corresponds to a vertex of degree $\le 2$ in $\Gamma$, then
contracting $\cD_i$ fiberwise yields a proper flat family of
surfaces $\phi:\bar \cV_{cont}\to S$ and a relative NC divisor
$\cD_{cont}$ on $\bar\cV_{cont}$ with a constant dual graph.

\item Assume that $\cD_i,\cD_j$ are irreducible components of $\cD$ such that
$\cD_i\cap\cD_j\cong S$ is a section of $\phi$.  Blowing up $\bar\cV$
along this section yields a flat family $\phi':\bar \cV'\to S$ together
with a relative NC divisor $\cD'=\phi'^{-1}(\cD)$ on $\bar\cV'$
with constant dual graph.
\enum
\elem

\bproof A proof of (a) can be found e.g.\ in \cite[Lemma
1.15]{FKZ3}. To deduce (b) it is enough to observe that locally
near $\cD_i\cap\cD_j$ our family is a product family. 
\eproof

We frequently use the following well known fact on the local triviality of  $\PP^1$-fibrations.

\blem\label{loctriv} 
Let $\phi:\cX\to S$ be a smooth morphism such that every fiber of $\phi$ is isomorphic to $\PP^1$. Then the following hold.

(a) $\phi$ is a locally trivial $\PP^1$-bundle in the \'etale topology. 

(b) If $\phi$ admits a section $\sigma:S\to\cX$ then $\phi$ is locally trivial in the Zariski-topology. 
\elem

\bproof
(b) is well known. It  can be shown e.g. along the same lines as Lemma 1.16 in \cite{FKZ3}. The assertion
(a) follows from  (b) since in a suitable \'etale neighborhood of a given  point $s\in S$ the map $\phi$ has sections.  
\eproof

Lemma \ref{simbd} yields the following result.

\blem\la{1.5} 
Let $\phi:\cV\to S$ be a family of normal 
surfaces, which admits an  NC completion
$(\bar\cV,\cD)$ with  constant dual graph $\G$. 
If  $\bdi \Gamma' &\rDotsto&\Gamma\edi$ is a birational transformation then locally in the \'etale topology there exists an NC completion $(\bar\cV' ,\cD')$ of $\phi$ such that $\cD'$ has constant dual graph $\Gamma'$.
\elem

\bproof
It is sufficient to treat the case where  $\Gamma'$ is obtained from $\Gamma$ by a single blowup or blowdown. The case of a single blowdown follows from Lemma \ref{simbd}(a), while in the case of an inner blowup part (b) of that Lemma is applicable. It remains to treat the case of an outer blowup in an irreducible component, say $\cD_i$ of $\cD$. Locally in the \'etale topology around a point $s\in S$ there are sections $S\to \cD_i$ not meeting the other components $\cD_j$ with $\cD_j\ne \cD_i$. Blowing up such a section yields the desired result. 
\eproof

We are now ready to give the proof of Theorem \ref{1.6}.

\bproof[Proof of Theorem \ref{1.6}]
By Theorem 2.15 in \cite{FKZ1} the dual graph $\Gamma$ of $(\bar\cV,\cD)$ can
be transformed into a standard graph by a sequence of
blowdowns and blowups. Applying Lemma \ref{1.5} the result follows. 
\eproof

Next we show that for completable families  $\phi:\cV\to S$ of affine ruled surfaces the morphism $\phi$ is affine. 

\bprop\label{fcontract} 
Let $\cV\to S$ be a completable family of normal surfaces   admitting an NC completion $(\bcV,\cD)$ with constant connected dual graph $\G$. 
\bnum[(a)]
\item If $\G$ has an extremal 0-vertex then   
there exists  a relative semi-ample divisor $A$ on $\bcV$ supported by $\cD$. 

\item The morphism $\cV\to S$ can be factorized 
into a proper relative simultaneous contraction $\cV\to \cV'$ (the Remmert reduction) and 
an affine morphism $\cV'\to S$. 
\enum  
\eprop

\bproof 
(b) follows immediately from (a). 

(a) If $\G$ has extremal 0-vertices then it is a tree. Letting $\cD=\sum_{i=0}^n\cD_i$, where $\cD_0=\cD_{01}$ corresponds to an extremal zero vertex of $\G$, we consider a $\QQ$-divisor 
$A=\sum_i a_i\cD_i$ supported by $\cD$, where $a_i>0$ $\forall i$. We let $a_0=1$, 
and then we 
choose the coefficients $a_i$ rapidly decreasing when the distance in $\G$ of $\cD_i$ to $\cD_0$ increases. Then in  each fiber, $A(s)\cdot\cD_i(s)>0$ for $i=0,\ldots, n$. 
\eproof

\bcor\la{semi-ample}
If $\phi:\cV\to S$ is a completable family of affine ruled surfaces then $\phi$ is an affine morphism. In particular if the base $S$ of the family is affine so is its total space $\cV$.
\ecor

\brem\label{simact} 
(1) Let $\phi:\cV\to S$ be a family of normal affine
surfaces which admits
an NC completion $(\bar\cV,\cD)$. If some fiber of $\phi$ is  an affine ruled surface then so is
every fiber of $\phi$. 
Indeed, the dual graph $\G$ of an NC completion of an affine surface $V$ can be transformed into a standard one $\G_{\rm std}$. Due to Remark \ref{cplusa}(1) the surface $V$ is affine ruled if and only if $\G_{\rm std}$ has extremal 0-vertices. Since the dual graphs $\G(s)$ of the fiber over $s\in S$ and then also $\G_{\rm std}(s)$ are constant the result follows. 

(2) Later on we will show that, under the assumptions
of (1), locally over $S$ the total space $\cV$ of the family carries a relative affine ruling. One can show that then it also carries locally a relative $\GG_a$-action.
\erem

\subsection{Extended divisors}\label{extdivi}
{\em In this subsection  $V$ denotes an affine ruled surface.}
In the sequel we use the following notation.

\bnota\label{resi} 
(a) Let $(\bV , D)$ be a semi-standard completion of $V$. Such a completion induces a $\PP^1$-fibration $\bar\pi:\bV\to \bB$, see Remark \ref{cplusa}(1).
In the case that $(\bV , D)$ is standard $\bar\pi$ was called in \cite{FKZ3} the {\em standard fibration} associated to $(\bV,D)$.  
Let $\pi:V\to B$ be the induced affine ruling. We note that by Lemma \ref{lemstd} any  affine ruling on $V$ appears in this way.

(b) The singularities of $\bV$ are all lying in the open part $V$. Thus, if $\rho:\tV\to\bV$ is the minimal resolution of singularities then $D$ can be considered as a divisor in $\tV$ denoted by the same letter $D$. We call $(\tV,D)$ a {\em resolved semi-standard completion} of $V$ and the induced $\PP^1$-ruling $\tilde\pi:=\bar\pi\circ\rho:\tV\to B$ again the associated standard ruling. 

(c) The horizontal section $C_1\cong \bar B$ of $\bar\pi$ (see Lemma \ref{lemstd}), considered as a curve in $\tV$, yields a horizontal section of $\tilde\pi$. 
Let $C_{0i}=\tilde\pi^{-1}(c_{0i})$, $i=1,\ldots, a$, denote the full fibers of $\tilde\pi$ contained in $D$, where 
$c_{0i}\in \bar B\setminus B$ are distinct points. They correspond to the extremal 0-vertices in $\G_D$ and are all adjacent to $C_1$. 
The vertex $C_1$ is connected to  further vertices, say
$C_{21}, \ldots, C_{2b}$.  
The latter represent fiber components
of the degenerate fibers, say, $\tilde\pi^{-1}(c_{21}),\ldots,
\tilde\pi^{-1}(c_{2b})$ of $\tilde\pi$. We let in the sequel
$$
C_0=\sum_{i=1}^a C_{0i}\and C_2=\sum_{j=1}^b C_{2j}.
$$
Since the component $C_{2i}$ meets the section $C_1$
transversally, it has multiplicity $1$ in the fiber
$\tilde\pi^{-1}(c_{2i})$. 
Hence  for every $i=1, \ldots, b$ the rest of
the fiber $\tilde\pi^{-1}(c_{2i})- C_{2i}$ can be blown down to a
smooth point. 
In this way we obtain a $\PP^1$-ruling (i.e.\ a locally trivial
$\PP^1$-fibration) $\psi: X\to \bar B$. Thus $\tV$ is obtained
from $X$ by a sequence of blowups  with centers on
$C_{2i}\backslash C_1$, $i=1, \ldots, b$, and at infinitesimally
near points. To simplify  notation, we keep the same letters
for the curves $C_{0i}$, $C_1$, and $C_{2i}$ on $\bV,\tV$ and for
their respective images in $X$.
\enota

In analogy with the case of Gizatullin surfaces \cite{FKZ2}, given a resolved standard completion $(\tV,D)$
of an affine ruled surface $V$, we associate to it the {\em extended divisor} $D_\ext$
and its {\em extended dual graph} $\G_\ext$ of  $(\tV,D)$ as follows.

\bdefi\label{exdi}
The reduced divisor
$$
D_{\ext}=D\cup \tilde\pi^{-1}( \{c_{21},\ldots, c_{2b}\})
$$
is called the {\em extended divisor} of $(\tV,D)$
and its dual graph $\Gamma_\ext$
the {\em extended graph}.
The connected components of $D_\ext- D$ are
called {\em feathers}, see \cite[2.3]{FKZ3}. 
\edefi

As shown there 
each feather of $D_\ext$ is a linear chain of smooth
rational curves on $\tV$
$$
\fF:\qquad \co{F_0}\lin \co{F_{1}}\lin\ldots\lin
\co{F_{ k}}\quad ,
$$
where the subchain $\fR=\fF- F_0=F_1+\ldots+ F_k$ (if non-empty)
contracts to a cyclic quotient singularity of $V$ and $F_0$ is
attached to some component $C$ of $D-C_0-C_1$. The
curve $F_0$ is called the {\em bridge curve}
of $\fF$.
For instance, an $A_k$-singularity on $V$ leads to
an $A_k$-feather, where  $\fR$ is a chain of $(-2)$-curves of
length $k$. In particular, $V$ can have at most cyclic quotient singularities. We note that this also follows from Miyanishi's Theorem \cite[Lemma 1.4.4]{Mi}.

In the case of a smooth surface $V$ one has
$k=0$ that is, every feather is an irreducible curve.
Otherwise the dual graphs of $\fR$ correspond
to the minimal resolutions of the singularities of $V$ and are independent of the
choice of a resolved completion. In contrast, the bridge
curve $F_0$, its neighbor in $D$, and its self-intersection
may depend on this choice.

We have the following uniqueness results for extended divisors.

\bprop\la{ext-cor} 
\bnum[(a)]
\item If the standard fibrations of
two standard completions $(\bV,D)$ and
$(\bV',D')$ of a Gizatullin surface $V$ induce the same  affine
ruling $V\to\AA^1$ then
there is an isomorphism of the associated extended divisors
$f:D_\ext\stackrel{\simeq}{\longrightarrow} D'_\ext$ and of their extended dual graphs
$\G_f:\G_\ext\stackrel{\simeq}{\longrightarrow} \G'_\ext$ sending
$D$ to $D'$ and $\G_D$ to $\G_{D'}$.

\item 
For a non-Gizatullin affine ruled surface $V$, 
up to an isomorphism the extended divisor 
$D_\ext$ and the extended graph $\Gamma_\ext$ 
do not depend on the choice of a standard 
completion $(\bV,D)$ of $V$.
\enum
\eprop

\bproof 
Performing an elementary transformation at an extremal zero
vertex $C_{0i}$ of $\G_D$ does not affect $D_\ext-C_{0i}$ while $C_{0i}$ is replaced by another smooth rational curve $C'_{0i}$. 
Thus under such an operation $D_\ext$ and $\G_\ext$ remain unchanged up to isomorphism. Now (a) and (b) follow from Lemma \ref{lemstd}(c). We note that (a) is also contained in Lemma 5.12 in \cite{FKZ3}.
\eproof

\bsit\label{simbl} In the setting of Notation \ref{resi}, the irreducible
components of $D_\ext-C_0-C_1-C_2$ in $\tV$
 are obtained via a sequence of blowups  starting
from the projective ruled surface $X$ as in  \ref{resi},
\be\la{blowupseq} 
\tV = X_m \to X_{m-1} \to \ldots \to X_2 \to
X_1=X\,, 
\ee 
with centers lying over $C_2
\backslash C_1$. Every component $F$ of $D_\ext-C_0-C_1-C_2$
is created by one of the blowups $X_{k+1}\to X_k$ in \eqref{blowupseq}. Since feathers do not contain $(-1)$-curves besides the bridge curves, the center of this blowup is necessarily a point lying on the image of $D$ in $X_k$. 
\esit

This justifies the following definition.

\bdefi\la{mother} (See \cite[2.3, 2.5]{FKZ3} or \cite[3.2.1]{FKZ4}.)
Suppose that a component $C$ of $D$ is mapped onto a curve $\bC$ in $X_{k-1}$. If under the blowup $X_k\to X_{k-1}$ the component $F$ of $D_\ext$ is created by a blowup on $\bC$ then $C$ is called a {\em mother component} of $F$. If the blowup takes place in the image of the point $p_F\in C$ then  $p_F$ is called the {\em base point} of $F$.
\edefi

The mother component $C$ of a component $F$ of a feather $\fF$ and its base point $p_F$ are uniquely
determined by $F$ (see \cite[Lemma 2.4(a)]{FKZ3}).
Indeed, otherwise $F$ would appear as the
exceptional curve of a blowup of an intersection point $C\cap C'$
in some $X_k$ with $C,C'\subseteq D$, so that $D$ looses
connectedness, which is impossible.

In contrast, a component $F=C$ of $D$ can have two mother components
(at most). For, if
$C\subseteq D$ is the result of an inner blowup of an intersection
point of two curves $C'$ and $C''$ of the image of $D$ in $X_k$,
then the proper transforms of
$C'$ and $C''$ in $D$ are both mother components of  $C$. The
curves $C_{0i}$, $C_{2j}$, and $C_1$ are orphans in that they
have no mother component. This distinction leads to the following definition.

\bdefi\label{tyco} Consider as before a resolved completion $(\tV,D)$ of an
affine ruled surface $V$.
We say that $C\subseteq D$ is a
{\em $*$-component}
if it has two mother components in $D$, and  a
{\em $+$-component}
otherwise. In particular, the curves $C_1$, $C_{0i}$, and $C_{2i}$ are
$+$-components.\edefi

\brem\label{-2g} 
Although the sequence (\ref{blowupseq}) is not unique, in general, we can replace it by a canonical one by blowing down on each step $X_k\to X_{k-1}$ simultaneously all
$(-1)$-curves in the fibers of the induced $\PP^1$-fibration $\pi_k:X_k\to B$ different from the components of the curve $C_2$. The latter contractions
are  defined correctly since, on each step,
no two $(-1)$-components in a degenerate fiber $\pi_k^{-1}(\pi_k(C_{2j}))$  different form
$C_{2j}$ are neighbors. This new sequence now only depends  on the completion $(\bV,D)$. 
\erem

\subsection{Extended divisors in families}\label{extdivifam}
Let now $\phi:\cV\to S$ be a family of normal affine surfaces. In
order to generalize the preceding construction to the relative
case, we will assume that $\phi$ is completable and admits a
simultaneous resolution of singularities i.e., a proper morphism
$\alpha: \cV'\to\cV$  such that the induced map $\phi':\cV'\to S$
is smooth and $\alpha(s): \cV'(s)\to\cV(s)$ is  a resolution of
singularities for every $s\in S$. Such a resolution will be called
{\em minimal}\/ if $\alpha(s): \cV'(s)\to\cV(s)$ is  a minimal
resolution of singularities for every $s\in S$.

In many cases, a simultaneuous resolution does not exist, see Remark \ref{extdefo} for a more thorough discussion. 
In the case where a simultaneous resolution does exist,
it can be chosen to be minimal due to the following  result.

\bprop\label{misire} 
Let $(\bcV, \cD)$ be an NC-completion of the family
$\phi:\cV\to S$. Assume that $\cV$ admits a simultaneous resolution
of singularities. Then locally in the \'etale topology of $S$ there is a minimal
simultaneous resolution of singularities $\alpha: \tcV\to\bcV$.
\eprop

\bproof Given a simultaneous resolution of singularities
$\cV'\to\cV$ there is an NC completion
$(\tcV,\cD)$ of $\cV'\to S$ by this same
relative NC divisor $\cD$ so that $\tcV\to\bcV$ is again a
simultaneous resolution of
singularities.

To make this resolution minimal let us fix a point $s_0\in S$. If
$\tcV(s_0)\to \bcV(s_0)$ is not a minimal resolution then there is
a $(-1)$-curve $E$ in the fiber $\tcV(s_0)$. By Lemma
\ref{deflem} near\footnote{In the \'etale topology.} $s_0$ there is a unique family of $(-1)$-curves
$\cE\to S$ in $\tcV$ with $\cE(s_0)=E$. Since $\cV$ is a family of affine surfaces, the fibers of $\cE$ are contracted in $\cV$ and then also in $\bcV$.
According to Lemma \ref{simbd}(a) contracting $\cE$ fiberwise
leads to a flat family $\tcV'$ over $S$. Replacing $\tcV$ by
$\tcV'$ and repeating the argument, after a finite number of steps
we arrive at a simultaneous resolution near
$s_0$ denoted again $\tcV\to \bcV$, which
is minimal in the fiber over $s_0$.

We claim that near $s_0$, the map $\tcV(s)\to \bcV(s)$ is a minimal resolution for
every $s\in S$. Indeed, consider the relative Hilbert scheme
$\HH_{\tcV/S}$. Every of its irreducible components is proper over
$S$ (see \cite{Gro0}), and the $(-1)$-curves in $\HH_{\tcV/S}$ form
a constructible subset, say, $A$. Then also the image of $A$ in
$S$ is constructible.

Hence, if in every neighborhood of $s_0$ there are points for
which the resolution $\tcV(s)\to \bcV(s)$ is not minimal, then there is a smooth
curve $T$, a morphism $\gamma:T\to S$, and a point $t_0\in T$ with
$\gamma(t_0)=s_0$ together with a $T$-flat family of curves
$\cC\subseteq \tcV\times_S T$ such that over $t\ne t_0$ the fiber
$C_t:=\cC(t)$ is a $(-1)$-curve in $\tcV(\gamma(t))$. Over $t_0$ we
obtain a curve $C= \cC(t_0)$ in $\tV=\tcV(s_0)$.

We claim that $C$ is contracted in $\bV$. To show this, let $p: \cC \to \bcV$ be the induced map. As $\cC$ is flat over $T$, the complement $\cC\backslash C$ is dense in $\cC$ and so by continuity $p(\cC\backslash C)$ is dense in $p(\cC)$.  
Since for $t\ne t_0$ the curve $C_t$ is contracted in $\bcV$ to a point,
the image $p(\cC\backslash C)$ and then also $p(\cC)$ are curves in $\bcV$. Thus $C$ is contracted to the point in $p(\cC)$ lying over $s_0$ in $\bV$, proving the claim. 

Because of $C.D=C_t.\cD(t)=0$ for $t\ne t_0$ the curve $C$ does not
meet $D$ and so is contained in $\tV\backslash D$. In
particular, since $\tV\to \bV$ is a minimal resolution of
singularities every
component $C_i$ of $C$ is smooth and $C_i^2\le -2$.

Letting $C=\sum n_iC_i$, $n_i>0$, and $K_t$ be the canonical
divisor on $\cV'(\gamma(t))$,  we get 
$$
-1=C_t.K_t=C.K=\sum n_iC_i.K\,,
$$
where $K=K_{t_0}$. Hence $C_i.K\le -1$ for some $i$.
This is only possible if $C_i^2\ge -1$, which contradicts the
minimality of resolution in the fiber over $s_0$.  Now the proof is completed.
\eproof

We need the following definition.

\bdefi\label{def-res-comp}
Let $\phi:\cV\to S$ be a flat family of normal surfaces,
$\phi':\tcV\to S$ a flat proper family of smooth surfaces and $\cD\subseteq \tcV$ a divisor. 
We call $(\tcV,\cD)$ a {\em resolved completion} of
$\phi:\cV\to S$ if the following conditions are satisfied.
\bnum[(a)]
\item 
$\phi':\cV':=\tcV\backslash \cD\to S$ is a minimal simultaneous
resolution  of singularities of $\phi:\cV\to S$;
\item 
$(\tcV, \cD)$ is a NC completion of $\phi':\cV'\to S$ as in Definition \ref{NCcomplet}; in particular it has a constant dual graph $\G$. 
\enum
We call the family $(\tcV,\cD)$ {\em standard or semi-standard} if so is $\G$. 
\edefi

In the next Proposition we show that one can organize the standard
morphisms of Notation \ref{resi} into a family.

\bprop\la{relstd}
Let $\phi:\cV\to S$ be a family of affine ruled surfaces,
which admits a resolved semi-standard completion $(\tcV,\cD)$ 
with associated map $\tilde\phi:\tcV\to S$.
Then there exists a factorization of $\tilde\phi$ as
\bdi
\tcV& \rTo^{\tilde\Pi} & \cB & \rTo &S\,,
\edi
where $\cB\to S$ is a family of
smooth complete curves and $\tilde\Pi$ induces in every fiber over $s\in S$ the standard
morphism $\tilde\pi_s$ from Notation \ref{resi}. 
\eprop

\bproof
The extremal zero vertices of $\G$ correspond to families of curves $\cC_{0i}\subseteq \cD$, $i=1,\ldots, a$. 
Their common neighbor $\cC_1\subseteq \cD$  is
a  family  of curves over $S$.
Given a point $s\in S$ we let $C_s$ and $\tV_s$ be the fibers of
$\cC=\cC_{01}$ and $\tcV$ over $s$, respectively. We consider a
semi-standard $\PP^1$-fibration $\tilde\pi_s:\tV_s\to \bB_s$ on $\tV_s$ as in \ref{resi}. The curve $C_s$ is the full fiber $\tilde\pi^{-1}(c_s)$ of $\tilde\pi_s$ over some point $c_s\in \bB_s$. The direct
image sheaf $R^1\tilde\pi_{s*}(\cO_{\tV_s}(mC_s))$  vanishes for $m\ge 0$, while $\tilde\pi_{s*}(\cO_{\tV_s}(mC_s))\cong \cO_{\bB_s}(m c_s)$.
Applying the Leray spectral sequence we obtain that
$H^1(\tV_s, \cO_{\tV_s}(mC_s)) \cong
H^1(\bB_s, \cO_{\bB_s}(m c_s))=0$ for $m\gg 0$ and all $s\in S$. 
Hence $R^1\tilde\phi_* (\cO_\tcV(m \cC))$ vanishes and is in particular locally free. By \cite[Chapt.\ III, Theorem 12.11]{Ha} the map
$$
\tilde\phi_*(\cO_\tcV(m \cC))_s\to H^0(\tV_s, \cO_{\tV_s}(m
C_s))
$$
is surjective for all $s\in S$. Shrinking $S$ we may
assume that there are sections
$$
\beta_0, \ldots ,\beta_N\in H^0(S, \tilde\phi_*(\cO_\tcV(m
\cC)))
$$ 
whose images in $H^0(\tV_s, \cO_{\tV_s}(m C))$
generate the latter vector space for every $s\in S$. Consider now
the morphism
$$
f=((\beta_0:\ldots :\beta_N), \tilde \phi): \tcV\to \PP^N\times S\,.
$$
By construction its restriction $f_s:\tV_s\to \PP^N$ to the fiber over $s$ factors through the standard morphism
\bdi
 \tV_s&\rTo^{\tilde\pi_s}& \bB_s&\rTo^{\gamma_s} \PP^N\,.
\edi
Here $\gamma_s$ denotes the map given by the linear system $|m\cdot c_s|$ on $B_s$. 

We claim that the sheaf $\cO_\cB=f_*(\cO_\tcV)$ is flat over $S$. As this is a local problem, we may suppose
in the proof of this claim 
that $S=\Spec A$ is affine. Using the vanishing of $H^1(\tV_s, \cO_{\tV_s})$ for all $s\in S$   by \cite[Chapt.\ III, Proposition 2.10]{Ha} the left exact functor 
$$
M\mapsto T(M)=H^0(\tcV, \cO_\tcV\otimes_AM)
$$
on finite $A$-modules is also right exact. 
Applying Proposition 12.5 in \cite[Chapt.\ III]{Ha} the natural map $T(A)\otimes_A M\to T(M)$ is an isomorphism. Thus the functor $M\mapsto T(A)\otimes_A M$ is exact and so $T(A)$ is a flat $A$-module. Consequently
$$
\cO_\cB:=f_*(\cO_\tcV)
$$
is a $S$-flat sheaf on $\PP^N\times S$, proving the claim.

The sheaf $\cO_\cB$  gives rise to a flat family of smooth curves $\cB\to S$ with fibers $\bB_s$. By
construction the morphism $\tilde\Pi=f:\tcV\to \cB$ induces over
each point $s\in S$ the standard fibration from Notation
\ref{resi}, proving the Lemma.
\eproof

\bdefi\label{stmofam} 
Let $\phi:\cV\to S$ be as before a family of affine 
ruled surfaces with a resolved semi-standard completion 
$(\tcV,\cD)$. The morphism $\tilde\Pi:\tcV\to \cB$ onto a
family of curves $\cB\to S$ constructed in the proof of Proposition
\ref{relstd} by means of the linear system $|m\cC|$, $m\gg 0$, will be called the {\em standard morphism} associated to $(\bcV,\cD)$. This is a $\PP^1$-fibration over $\cB$ which induces in each fiber over
$s\in S$ the standard $\PP^1$-fibration. 
\edefi

\bsit\label{order.2} The vertex of $\G$ which corresponds to
$\cC_1$ has neighbors given by $\cC_{01}, \ldots, \cC_{0a}$ and
further ones given (locally in $S$) by smooth  families of
rational curves $\cC_{21},\ldots, \cC_{2b}\subseteq\tcV$ over $S$.
The families $\cC_{01},\ldots, \cC_{0a}$  arise as preimages under
$\tilde\Pi$ of sections, say, $\gamma_{0i}:S\to\cB$,
$i=1,\ldots,a$, while $\cC_{2j}$ are projected in $\cB$ to
sections, say, $\gamma_{2j}:S\to\cB$, $j=1,\ldots, b$. Moreover
$\cC_1\cong \cB$ is a section of the standard morphism $\tilde\Pi$.
On the curve $\bar B_s=\cB(s)$ the sections $\gamma_{0i}$ and
$\gamma_{2j}$ yield the points $c_{0i}=\gamma_{0i}(s)$ and
$c_{2j}=\gamma_{2j}(s)$ (cf.\ Definition \ref{resi}). 
\esit

The following Factorization Lemma will be useful in the sequel.

\blem \la{factlem} {\bf (Factorization Lemma)} Let $\cV\to S$ be a
family of affine ruled surfaces, which admits a resolved standard
completion $(\tcV, \cD)$
(see Definition \ref{def-res-comp}).
Then locally in the \'etale topology of $S$ there is a
factorization of the associated standard morphism $\tilde \Pi:\tcV\to\tB$
as
\be\la{fact} 
\tcV=\cX_m\to \cX_{m-1}\to \ldots \to
\cX_1\to \cB\,, 
\ee 
where 
\bnum[(a)]
\item every morphism $\cX_i\to\cX_{i-1}$ is
a blowup of a section $\gamma_i:S\to \cX_{i-1}$ contained in the image of $\cD$ in $\cX_ {i-1}$ with  the exceptional divisor
$\cE_i\subseteq \cX_i$ being a  relative (-1)-curve over $S$,

\item $\cB\to S$ is a smooth family of
complete curves over $S$,  and

\item $\cX_1\to \cB$ is a locally trivial
$\PP^1$-bundle. 
\enum
\elem

\bproof 
On the fiber over some point
$s$ we consider the standard fibration
$$
\tilde\pi_s=:\tilde\Pi|\tV_s: \tV_s=\tcV(s)\to  \bB_s \,.
$$
We let $C_{0i,s}=\cC_{0i}(s)$, $C_{1,s}=\cC_{1}(s)$, and
$C_{2j,s}=\cC_{2j}(s)$ be the respective fibers over $s$ so that $C_{2j,s}$ is contained in a fiber $\tilde\pi_s^{-1}(c_{2j,s})$ over some point $c_{2j,s}$ of $\bB_s$.
Blowing down the divisors $\tilde\pi_s^{-1}(c_{2j,s})-C_{2j,s}$, $j=1,
\ldots, b$, we arrive at a locally trivial $\PP^1$-bundle, say,
$\Psi: X\to \bar B$ such that $\bV$ is obtained by a sequence of
blowups  of $X$ as in (\ref{blowupseq}) with centers on
$\bigcup_{j=1}^b\Psi^{-1}(c_{2j})\backslash C_1$ and
infinitesimally near points.\footnote{As in Notation
\ref{resi} we  use  the same letters for $C_{0i}$, $C_1$
and for their images in $X$.}

In particular, in the last blowup $X_m=\tV\to X_{m-1}$ in
(\ref{blowupseq}) there is a $(-1)$-curve, say, $E_m\subseteq \tV$
which is blown down in $X_{m-1}$. By Lemma \ref{deflem}(a) near $s_0$ there is a
family of $(-1)$-curves $\cE_m\subseteq \tcV$ inducing $E_m$ over
$s_0$. Applying Lemma \ref{simbd}(a) we can blow down $\cE_m$ and
obtain a morphism $\cX_m=\tcV\to \cX_{m-1}$ inducing $X_m\to
X_{m-1}$ in the fiber over $s_0$. Repeating this procedure we
arrive at a factorization (\ref{fact}).

It remains to prove that $\cX_1\to\cB$ is a locally trivial
$\PP^1$-fibration. Let us show first that the morphism
$\cX_1\to\cB$ is flat. Indeed, using in every step of our
construction Lemma \ref{simbd}(a)  $\cX_i\to S$ is a flat
morphism. In particular $\cX_1\to S$ is flat. Since for every
$s\in S$ also $\cX_1(s)\to \cB(s)$ is flat, the flatness of
$\cX_1\to\cB$ follows from \cite[Corollary 6.9]{Ei}.

Now the fact that $\cX_1\to\cB$ is a locally trivial
$\PP^1$-fibration is a consequence of Lemma \ref{loctriv}. 
For $\cC_1\to \cB$ is an isomorphism and so the inclusion
$\cC_1\hookrightarrow \cX_1$
yields a section of $\cX_1\to\cB$. 
\eproof

By construction all blowups in the sequence (\ref{fact}) take
place over disjoint sections $\gamma_{21}(S),\ldots,\gamma_{2b}(S)$
of $\cC_2\setminus\cC_1=\bigcup_{j=1}^b
\cC_{2j}\backslash \cC_1$  over $S$ in $\cX_1$. In analogy with the absolute
case
we introduce the extended divisor of $(\bar\cV,\cD)$ as
\be\la{relext}
\cD_\ext= \cD\cup \tilde\Pi^{-1} \left( \gamma_{21}(S)\cup \ldots
\cup \gamma_{2b}(S)\right)\,
\ee
(see Definition \ref{exdi}) .
We emphasize that this is not, in general, a relative SNC
divisor in the sense of Definition \ref{reldiv} (see Example \ref{jump} below).
However, we have the following result.

\bcor\label{unnamed} 
Under the assumptions of Lemma \ref{factlem} 
the extended divisor $\cD_\ext$ is flat over $S$, and each fiber
$\cD_\ext(s)$ over $s\in S$ is just the extended divisor of $(\tV_s, D_s)$. 
\ecor

\bproof 
Locally $\cD_\ext$ is the set of zeros of a non-zero divisor on $\tcV$ that is also a non-zero divisor in each fiber. Thus the first part follows e.g.\ from \cite[Corollary 6.9]{Ei}. 
The second part follows from the fact that blowing up a section of $\tcV\to S$ commutes with taking the fiber over $s$.
\eproof

\brem\la{types-blowups}  In analogy with the absolute case we call a component  $\cA$ of the extended divisor $\cD_\ext$ from \eqref{relext} a {\em boundary component} if it is in $\tcV$  a component of $\cD$, and otherwise a {\em feather component}.
Since the dual graph of $\cD(s)$ is constant,  $\cA$ is a boundary or a feather component if and only if its fiber $\cA(s)$ over some point $s\in S$ is. However the neighbouring components of $\cA(s)$ may change in nearby fibers. We call this phenomenon ``jumping''\ of feathers.
\erem

Let us give an example where this actually happens.

\bexa\la{jump}  
Letting $V$ be the Danilov-Gizatullin surface with dual zigzag $\G=[[0,0,-2,-2]]$ we consider the trivial family $\cV=V\times S$ over $S=\AA^1$. 
We construct a family of completions $(\bar\cV, \cD)$ over $S$ with constant dual zigzag $\G$, where the dual graph $\cD_\ext(s)$ of the extended divisor is not constant on $S$. Actually we will see that the dual graphs of $\cD_\ext(s)$, $s\ne0$, and $\cD_\ext(0)$ are as follows:
\vspace{6mm}
\be\la{jumpgraph}
\cD_\ext(s):
\quad \cou{0}{C_0}\lin\cou{0}{C_1}\lin
\cou{\qquad-2}{C_2}\nlin\cshiftup{F_1}{-1}\lin\cou{\,\quad-2}{C_3}\lin\cou{-1}{F_2}\quad, \quad
\qquad\cD_\ext(0):
\quad \cou{0}{C_0}\lin\cou{0}{C_1}\lin
\cou{-2}{C_2}\lin\cou{\qquad-2}{C_3}\nlin\cshiftup{F_1}{-2}\llin\cou{-1}{F_2}
\ee
The construction starts  from the quadric  $X_1=Q=\PP^1\times\PP^1$  with the $\PP^1$-fibration given by the first projection $X_1\to B=\PP^1$ and with the curves
$$
C_0=\{\infty\}\times\PP^1,\qquad C_1=\PP^1 \times\{\infty\},\quad\text{and}\quad
C_2=\{0\}\times\PP^1\,.
$$
Blowing up $X_1$ at the point $(0,0)$ creates a feather $F_1$ on the blown up surface $X_2$.
Letting $\cX_2=X_2\times\AA^1$, $\cC_i= C_i\times\AA^1$, $i=0,1,2$, and  $\cF_1= F_1\times\AA^1$, we blow up $\cX_2$ along a ``diagonal'' section $\gamma_3:S\to \cC_2$,
which meets $\cF_1$ over $s=0$ only. By this blowup in every fiber over $s$ the boundary component $C_3$ is created. 
Blowing up a  section $\gamma_4(S)\subseteq \cC_3$, which does not meet $\cC_2\cup\cF_1$, we obtain a second feather $\cF_2$ of $\cD_\ext$ on the new threefold $\bcV=\cX_4$.
Thus the feather $\cF_1(s)$, $s\ne 0$, jumps from $C_2$ to $C_3$
over the point $s=0$ of the base as indicated by the dual graphs. \eexa

We end this section with a remark on the existence of a resolved completion in our setup.

\brem\la{extdefo}
(1)
Let $V$ be a normal affine ruled surface with a 
standard completion $(\bV, D)$ and associated 
standard fibration $\bV\to \bar B$.
Then every deformation $\cV\to S$ of $V$ over an 
Artinian germ $S$ can be extended to a deformation 
$(\bcV,\cD)\to S$ of the completion  $(\bV,D)$.
Indeed, the infinitesimal deformations and obstructions of
$(\bV,D)$ are given by
$$
T^1=H^1(\bV, R\cHom_\bV(\Omega^1_\bV\langle D\rangle, \cO_\bV))
\and T^2=H^2(\bV, R\cHom_\bV(\Omega^1_\bV\langle D\rangle, \cO_\bV))
$$
respectively, see \cite{Se}. Consider the triangle
in the derived category
\be\la{eseq}
0\to \Theta_\bV\langle D\rangle
\to R\cHom_\bV(\Omega^1_\bV\langle D\rangle, \cO_\bV)\to \cG^\bullet\to 0\,.
\ee
The first and the second cohomology of $\cG^\bullet$ are just
$$
T^1_\loc=\Ext^1_V(\Omega^1_V, \cO_V))
\and T_\loc^2=\Ext^2_V(\Omega^1_V, \cO_V))\,,
$$
which control the deformations and obstructions of the affine part,
respectively (cf.\ \cite{Se}).
Using the long exact cohomology sequence of (\ref{eseq})
we get an exact sequence
$$
T^1\to T^1_\loc\to H^2(\bV, \Theta_\bV\langle D\rangle)
\to T^2\to T^2_\loc\,.
$$
The term in the middle is dual to 
$\Hom_\bV(\Theta_\bV\langle D\rangle, \omega_\bV)$
and so vanishes as the restriction of an element of the latter group to a generic fiber of $\bV\to \bar B$ vanishes.

Thus by standard deformation theory (see\cite{Se})
the functor assigning to a deformation of $(\bV,D)$ the deformations of the open part is formally smooth.

(2) Consider, for instance, a surface $V$ which has a cyclic quotient singularity
with several components in the versal deformation. Its  versal
family $\cV\to S$ admits at least formally a simultaneous completion by (1) but does
not admit a simultaneous resolution. See \cite{KSB} for examples of cyclic quotients for which the versal deformation space is not
irreducible and does not coincide with the Artin component.
\erem

\section{Deformation equivalence}\label{defequiv}

\subsection{Normalized extended graph and
deformation equivalence}

The purpose of this and the next sections is to characterize in
combinatorial terms the deformation equivalence for affine ruled surfaces, which we introduce as follows.

\bdefi\label{deeq}
We say that two normal surfaces $V$ and $V'$ are {\em deformations of each other}
if there exists a flat family of surfaces $\phi:\cV\to S$ over a connected
base $S$ with the following properties.
\bnum
\item[(a)]  $\phi$ admits a resolved completion
(see Definition \ref{def-res-comp}), and
\item[(b)]
$V\cong \cV(s)$ and
$V'\cong \cV(s')$ for some points $s,s'\in S$. 
\enum
This
generates an equivalence relation called {\em deformation
equivalence}: $V\sim V'$ if and only if there exists a chain
$V=V_1,\ldots,V_n=V'$ such that $V_i$ and $V_{i+1}$ are
deformations of each other for every $i=1,\ldots,n-1$.
\edefi

In the sequel $V$ will be a normal affine ruled surface. 
Let $(\bV,D)$ be a standard completion of $V$,
and  let $\tV\to\bV$ be the minimal
resolution of singularities.  The extended divisor $D_\ext$ as in Definition \ref{exdi} and the
extended graph $\Gamma_\ext=\G(D_\ext)$ depend in general 
on the choice of a completion $(\bV,D)$. We
associate to $\Gamma_\ext$ now another graph $N(\Gamma_\ext)$ called
the {\em normalized extended graph} which will turn out to be independent of the choice of completion and thus is an invariant of the affine surface.

\bdefi\label{normextgr} 
Given a component $C$ of $\G$ we let $\delta_C$ denote  the number of feather components $F$ of $\Gamma_\ext$ with mother
component $C$ (see Definitions \ref{exdi} and \ref{mother}). 
The {\em normalized extended graph}\/
$\Delta=N(\Gamma_\ext)$ of $(\bV, D)$ is the weighted graph
obtained from $\Gamma=\G_D$ by attaching to every component $C$ of $\Gamma$ exactly $\delta_C$ extremal $(-1)$-vertices called the \em{feathers of} $N(\Gamma_\ext)$.
\edefi

Thus $\Gamma_\ext$ and
$N(\Gamma_\ext)$ contain both $\Gamma$ as a distinguished subgraph
and have the same number of feather components, and even the same
number of them with a given mother component.  Furthermore,
every feather of $N(\Gamma_\ext)$ consists of
a single extremal $(-1)$-vertex, and these vertices
are in one-to-one correspondence with
the  feather components of $\Gamma_\ext$. 
For instance,
for the family $(\cX(s),\cD(s))$ as in Example \ref{jump} we have $\delta_{C_2}=\delta_{C_3}=1$ at any point $s\in S$, so the normalized graph $\Delta$ is in both cases the dual graph on the left  in \eqref{jumpgraph}.

\brem\la{ext-next} The normalized extended graph $N(\Gamma_\ext)$ can
be uniquely recovered from the extended graph $\Gamma_\ext$ via the simultaneous contraction
procedure described in \ref{simbl}.
\erem

In the case of a Gizatullin surface $V$ with a standard 
completion $(\bar V,D)$ the boundary zigzag $\G$ can be reversed
by moving the two zeros to
the other end as in (\ref{reverse}). Since $\G$ is contained in the extended divisor $\Gamma_\ext$ and also in $\Delta=N(\Gamma_\ext)$, we can perform the same operation in $\Gamma_\ext$ and in $\Delta$. In this way we obtain a new normalized
extended graph $\Delta^\vee=N(\Gamma_\ext^\vee)$
called the
{\em reversion} of $\Delta$.

\bexa\label{dgsu} Recall that every Danilov-Gizatullin surface $V=V_n$
(see the Introduction)
can be realized as
the complement of an ample section $\sigma$ in a
Hirzebruch surface $\Sigma_d$ with self-intersection
$\sigma^2=n>d$. By a theorem of Danilov-Gizatullin
\cite[Theorem 5.8.1]{DaGi}\footnote{See
also \cite{CNR, FKZ5}.} the
isomorphism class of $V_n$ depends only on $n$ and not on the
choice of $\sigma$ or of the concrete Hirzebruch surface
$\Sigma_d$. According to Example 1.22 in \cite{FKZ3}, for every
$r$ in the range $1,\ldots,n-1$ the surface $V_n$ admits a standard completion $(\bV_{n,r}, D)$ with extended
graph

\be\label{DGext} \G_{\ext,r}:\qquad \cou{C_0}{0}\lin\cou{C_1}{0}\lin
\cou{C_2}{-2}\lin\ldots\lin\cou{\quad\qquad C_{r+1}}{-2}
\nlin\cshiftup{\fF_1}{-r}
\vlin{20}\cou{C_{r+2}}{-2}\lin\ldots\lin\cou{\quad\qquad
C_{n}}{-2} \nlin\cshiftup{\fF_0}{-1}\qquad\quad,
\ee
where the bottom line $\G_D$ corresponds to the
boundary zigzag $D$, the feather $\fF_1$ consists of a single
$(-r)$-component $F_1$, and $\fF_0$
of a single $(-1)$-component (for
$r+1=n$ both feathers $\fF_0$ and $\fF_1$ are attached to
$C_{n}$). For every $r=1,\ldots,n$ the mother component of $F_1$
is $C_2$. Thus all these completions have the same normalized extended graph
$N(\G_{\ext,r})=\G_{\ext,1}$. In particular, $\G_{\ext,1}$ is as well
the  normalized extended graph
of the reversed standard completion
$(\bV_{n,r}^\vee, D^\vee)$ with the extended graph $\G_{\ext,n-r+2}$.
\eexa

In the next result we show that $N(\G_\ext)$ is essentially independent of the choice of the completion $(\bV,D)$ of an affine ruled surface $V$ and thus is an invariant of $V$.

\bthm\la{extinv} 
Let $V$ be a normal affine ruled surface.
\bnum[(a)]
\item If $V$ is not Gizatullin then up to an isomorphism
its normalized extended graph $N(\Gamma_\ext)$
is independent of the choice
of a resolved standard completion $(\tV,D)$ of $V$.\footnote{Hence
this graph $N(\Gamma_\ext)$ is a combinatorial  invariant of the affine ruled surface $V$.}

\item[(b)] If $V$ is Gizatullin then its normalized extended graph
$N(\Gamma_\ext)$ is uniquely determined by $V$ up to reversion.
\enum
\ethm 

\bproof
Assertion (a) follows from Corollary \ref{ext-cor}(b). Indeed, due to this corollary for a non-Gizatullin affine ruled surface $V$ the extended graph $\G_\ext$,
and hence also the normalized extended graph $N(\G_\ext)$, does not depend on the choice of a standard completion.

In the smooth case the proof of (b) is a consequence of \cite[Corollary 3.4.3]{FKZ3}. For the normal case we provide a proof in Theorem \ref{matchthm} below. 
\eproof

Now we come to the main theorem of our paper. 
In terms of the normalized extended graphs deformation equivalence
can be characterized as follows.

\bthm\la{main1} 
Two  affine ruled surfaces $V$ and $V'$ with resolved standard completions $(\bV,D)$ and $\bV',D')$ are deformation equivalent if and only if the following two conditions hold:
\bnum[(i)]
\item The associated normalized extended graphs are isomorphic  or, in the case of Gizatullin surfaces, are isomorphic 
up to reversion.
\item  The horizontal curves $C_1\subseteq D$ and $C_1'\subseteq D'$ have the same
genus.
\enum
\ethm

We establish the `only if'-part in subsection \ref{definv}, and  the
`if'-part in Corollary \ref{tercor} in section \ref{sec-versal}.

\subsection{A general matching principle}\label{gmp}
Part (b) of Theorem \ref{extinv} follows immediately from
Theorem \ref{matchthm}(b) below. To formulate this theorem we introduce the following notation.

\bnota\la{star-plus} 
Let $C$ be a component of $D-C_0-C_1$.
Then $C_1$ is sitting on a unique branch $\G_{D'}$ of the tree $\G_D$ at the vertex $C$; note that $D'$ contains a mother component of $C$.
We let $C\cap D'=\{\infty_C\}$.
For a $+$-component $C$ (see Definition \ref{tyco}) we let
$C^*=C\backslash \{\infty_C\}=C\backslash D'$.

A $*$-component $C$ has two mother components,
say, $C'\subseteq D'$ and $C''\subseteq D-D'$.
For such a component $C$ we let
$C^*=C\backslash (D'\cup D'')=C\setminus\{\infty_C,0_C\}$, where $D''\neq D'$ is the branch of $D$ at $C$ containing $C''$ and $\{0_C\}=C\cap D''$. 
\enota

Next we define the {\em configuration invariant} of $V$. In the particular case of smooth
Gizatullin surfaces, this invariant
was introduced in \cite[\S 3.2]{FKZ4}.

\bdefi\la{config-sp} (cf.\ \cite[(3.2.3)]{FKZ4}) Given a component $C$ of $D-C_0-C_1$
we  consider the cycle on $C$
$$
Q_C=\sum_F p_F\,,
$$
where $F$ runs over all feather components with
mother component $C$, and $p_F\in C$
is the base point of $F$ (see
Definition \ref{mother}).
From the fact that every feather has just one mother component it follows that
the cycle $Q_C$ is contained in $C^*$.
Letting $\delta_C=\deg Q_C$ be the number of feathers with mother component $C$,  we consider
the Hilbert scheme  $\HH_{\delta_C}(C^*)$ of subschemes
of length $\delta_C$ in $C^*$, or, equivalently,
of effective zero cycles of degree $\delta_C$ on $C^*$. This space is the quotient of
the Cartesian power $(C^*)^{\delta_C}$ modulo the symmetric group $\cS_{\delta_C}$.
Let $\Aut(C^*, \infty_C)$  denote the group of automorphisms
of the curve $C\cong \PP^1$  fixing the point $\infty_C=C\cap D'$ and, if $C$ is a $*$-component,
also the point $0_C=C\cap D''$ (see \ref{star-plus}).
This group acts diagonally on $(C^*)^{\delta_C}$
commuting with the $\cS_{\delta_C}$-action.
Hence it also acts on the Hilbert scheme $ \HH_{\delta_C}(C^*)$. Let us consider the quotient
$$
\fC_{\delta_C}(C^*)= \HH_{\delta_C}(C^*)/\Aut(C^*,\infty_C)\,.
$$
The cycle $Q_C$
defines a point denoted by the same letter in the space
$\fC_{\delta_C}(C^*)$.
Letting $$ \fC(\bV,D)=\prod_{C\subseteq D-C_0-C_1}\fC_{\delta_C}(C^*)
\,$$ we obtain a point
$$
Q(\bV,D)=\{Q_C\}_{C\subseteq D-C_0-C_1}\in \fC(\bV,D)
$$
called the {\em configuration invariant} of $V$.
\edefi

\brem\la{extra}
Performing in $(\tV,D)$ elementary transformations 
with centers at the components $C_{0i}$ of $C_0$ we neither
change $\tilde\Pi$ nor the extended divisor 
(except for the self-intersection index $C_1^2$) and 
thus leave the data $\delta_C$ and $Q(\bV,D)$ invariant. 
Hence we can define these invariants also for any 
semi-standard completion $(\bV, D)$ of $V$ by
sending the latter via elementary transformations with centers on $C_0$ into a standard completion. 
\erem

In order to show that the configuration
invariant is independent of the choice of a standard completion we need the
following definition.

\bdefi\label{fco} Let  $V$ be an affine ruled surface.
Given two standard completions $(\bV, D)$ and
$(\bV',D')$ of $V$, we consider the birational map $f:\bV\to\bV'$,
which extends the identity map of $V$. We distinguish between the
following two cases.

(1) If $V$ is non-Gizatullin then according to Proposition \ref{ext-cor}(b) the extended divisor of a standard completion is uniquely determined. In other words, 
$f$ induces a canonical isomorphism of extended divisors $D_\ext$ and $D'_\ext$.
The component in $D'$ corresponding under this isomorphism to the component $C\subseteq D$
will be denoted  $C^f$.

(2) Assume further that $V$ is a Gizatullin surface,
and let $D=C_0\cup
\ldots\cup C_n$ and $D'=C'_0\cup \ldots\cup C'_n$
be the standard
zigzags of the corresponding completions.
Then $(\bV', D')$ is
symmetrically linked either to the completion $(\bV,D)$,
or to its reversion
$(\bV^\vee, D^\vee)$, see \cite[2.2.1-2.2.2]{FKZ4}.\footnote{The latter means that there is a sequence of blowups and blowdowns transforming $(\bV', D')$ into one of $(\bV,D)$ or
$(\bV^\vee, D^\vee)$ and inducing an isomorphism of the corresponding dual graphs, which can be written symmetrically as $(\gamma,\gamma^{-1})$, where $\G_{D'}\to\G'$ is a birational transformation.}
In the first case we set $C_i^f=C'_i$ and in
the second one $C_i^f=C'_{i^\vee}$, where
$$
i^\vee:=n-i+2\,.
$$
If $(\bV', D')$ is symmetrically linked to both
completions $(\bV,D)$ and $(\bV^\vee, D^\vee)$, then we
define $C_i^f$ to be either $C'_i$ or $C'_{i^\vee}$.
\edefi

The following theorem in the particular case of smooth Gizatullin surfaces
is proven in \cite[Proposition 3.3.1]{FKZ4}.

\bthm\label{matchthm} {\bf (General Matching Principle)} Let
$V$ be an affine ruled surface with two standard completions $(\bV, D)$
and $(\bV',D')$, and let $C$ be a component of $D$ with the corresponding component
$C^f$ of $D'$. Then the following hold.
\bnum[(a)]
\item $C$ is a $*$-component if and only if $C^f$ is;

\item there is a canonical isomorphism $C\cong C^f$ 
mapping $C^*$ to $(C^f)^*$ and $\infty_C$ to $\infty_{C^f}$.
Under this isomorphism the cycle $Q_C$ on $C$
is mapped onto the cycle $Q_{C^f}$ on $C^f$;

\item the induced isomorphism
$\fC(\bV,D)\cong \fC(\bV',D') $
sends the configuration invariant $Q(\bV,D)$ to $Q(\bV',D')$.
\enum
\ethm

\bproof
Clearly (c) is a consequence of (a) and (b).
If $V$ is not a Gizatullin surface, then
the extended divisors $D_\ext$ and
$D'_\ext$ are  isomorphic, see Proposition \ref{ext-cor}(b).
Hence (a) follows   in this case.
Since the cycles $Q_C$ can be
read off from this extended divisor,
also (b) follows.

In the case of a Gizatullin surface $V$, the birational map
$\bdi(\bV,D)&\rDotsto &(\bV',D')\edi$ induced by $\rm id_V$ can be uniquely decomposed
into a sequence
$$
D=Z_1\sto Z_2\sto  \ldots\sto Z_n=D'\,,
$$
where each $Z_i$ is a semi-standard zigzag and each step is either
\bnum[(i)]
\item the reversion of a standard zigzag, or

\item an elementary transformation at an extremal $0$-vertex,
\enum
see \cite[Theorem 1]{DaGi} or \cite[Proposition 2.3.3]{FKZ4} for the existence part.
A transformation of type (ii) does not alter the extended divisor except for the weight  $C_1^2$, so it preserves the configuration invariant.
Hence we are done in this case as before. It remains to deduce (a) and (b)
in the case, where $(\bV',D')$ is the reversion of $(\bV,D)$.
This is the content of the following proposition.
\eproof

\bprop\la{matchprop}
Let $V$ be a Gizatullin surface with standard completion
$(\bV, D)$ and standard zigzag $D=C_0\cup \ldots\cup C_n$. Let further
$(\bV^\vee,
D^\vee)$ be the reversed completion with
the reversed standard zigzag $D^\vee=C^\vee_0\cup \ldots\cup C^\vee_n$.
Then the
following hold.
\bnum[(a)]
\item $C_i$ is a $*$-component if and only if $C^\vee_{i^\vee}$ is;
\item there is a natural isomorphism $C_i\cong C^\vee_{i^\vee}$ mapping $C^*$ to  $(C^\vee)^*$ and
$\infty_C$ to $\infty_{C^\vee_{i^\vee}}$. Under this isomorphism the
cycle
$Q_{C_i}$ is mapped onto the cycle $Q_{C^\vee_{i^\vee}}$.
\enum
\eprop

The proof is given in \ref{1.91} below.
Let us recall first the
 following notation. 

\bnota
Given a Gizatullin surface $V$ with boundary zigzag
$D=C_0+C_1+\ldots+C_n$, for $t\in\{1,\ldots,n\}$
we let $D_\ext^{\ge t}$ denote the branch of $D_\ext$ at $C_{t-1}$ containing $C_t$.
Moreover we let $D^{\ge t}=D\cap D_\ext^{\ge t}$ and  $D_\ext^{> t}=D_\ext^{\ge t}-C_t$ (see \cite[3.2.1]{FKZ4}).
\enota

The proof of Proposition \ref{matchprop} is similar to the proof for smooth  Gizatullin surfaces given in  \cite[Prop.\ 3.3.1]{FKZ4}. As in {\em
loc.cit.} the main tool is the {\em correspondence fibration}. Let us
recall this notion from  \cite[3.3.2-3.3.3]{FKZ4}.

\bdefi\la{corrfib} Using inner elementary transformations we can move the pair of
zeros in the zigzag $D=[[0,0,w_2,\ldots, w_n]]$ several places to
the right. In this way we obtain a new resolved completion
$(W,D_W)$ of $V$ with boundary zigzag $$D_W=[[w_2,\ldots, w_{t-1}, 0, 0, w_t,
\ldots,w_n]]$$ for some $t\in\{2,\ldots,n+1\}$. For $t=2$,
$D_W=D=[[0, 0,w_2,\ldots, w_n]]$ is the original zigzag, while for
$t=n+1$, $D_W=D^\vee=[[ w_2,\ldots, w_n,0, 0]]$ is the reversed
one. The
new zigzag $D_W$ can also be written as
$$
D_W=C_n^\vee \cup \ldots\cup C_{t^\vee}^\vee \cup C_{t-1} \cup
C_t\cup \ldots \cup C_n=D^{\ge t-1}\cup D^{\vee \ge t^\vee}\,,
$$
where for all $t-1\le i\le n$ and
$t^\vee\le j\le n$ we identify $C_i\subseteq \tV$ and $C_j^\vee\subseteq
\tV^\vee$ with their proper transforms in $W$.
In $W$ we have $C_{t-1}^2=C^{\vee 2}_{t^\vee}=0$.
Likewise in \cite[3.3.3]{FKZ4} the map
$$
\psi: W\to \PP^1
$$
defined by the linear system $|C_{t-1}|$ on $W$  will be called
the  {\em correspondence fibration} for the pair of curves $(C_t,
C_{t^\vee}^\vee)$.
The components $C_t$ and $C_{t^\vee}^\vee$ represent sections of
$\psi$. Their projections to the base  $\PP^1$ yield isomorphisms
$C_t\cong C_{t^\vee}^\vee\cong \PP^1$; in what follows
we identify points of all three curves under these isomorphisms.

Since the feathers of $D_\ext$ and $D_\ext^\vee$ are not
contained in the boundary zigzags they are not contracted in $W$.
We denote their proper transforms in $W$ by the same letters. It
will be clear from the context where they are considered.
We observe that every such feather $\fF$ can be written as
$\fF=F+\fR$, where $F$ is a unique component of $\fF$ with $F\cdot D_W\ge 1$
and  $\fR$  is the exceptional divisor of the minimal resolution in
$W$ of a cyclic quotient singularity of $V$ with $\fR\cdot D_W=0$.  \edefi

The following lemma is proven in \cite[Lemma 3.3.4]{FKZ4} for smooth
Gizatullin surfaces; the proof in the general case is similar and can be left to the reader.

\blem\label{mem} Under the assumptions of Proposition \ref{matchprop},
with the notation as in \ref{corrfib} the following hold.
\bnum[(a)]\item The divisor $D_\ext^{\ge t+1}$ is
contained in some fiber $\psi^{-1}(q)$,  $q\in\PP^1$. Similarly,
$D_\ext^{\vee \ge t^\vee +1}$ is contained in some fiber
$\psi^{-1}(q^\vee)$. The points $q$ and $q^\vee$ are uniquely
determined unless $D_\ext^{\ge t+1}$ and $D_\ext^{\vee\ge t^\vee
+1}$ are empty, respectively. \item A fiber $\psi^{-1}(p)$ contains
at most one component $C$ not belonging to $D_\ext^{> t}\cup
D_\ext^{\vee> t^\vee}$. Such a component $C$ meets both $D^{\ge
t}$ and $D^{\vee\ge t^\vee}$.\enum \elem

The next lemma is crucial in the proof of Proposition
\ref{matchthm}; see \cite[3.3.6 and 3.3.9]{FKZ4}
for the case of a smooth Gizatullin surface. The proof in the general case is far more involved.

\blem\label{matchlem} Given a point $q\in \PP^1$, we let $F_0,\ldots, F_k$
denote the feather components of the extended divisor $D_\ext$   with mother component $C_t$ contained
in the fiber $\psi^{-1}(q)$. Assuming that there exists at least one such component,
i.e. $k\ge 0$, the following hold. 
\bnum[(a)]

\item  With a suitable enumeration, the components $F_0,\ldots, F_k$ form a chain in
$D^{\ge t}_\ext$ contained in some feather $\fF$ of $D_\ext^{\ge t}$.

\item
Let $F_0$ have the smallest distance to $C_t$ in the
chain above. Then the fiber $\psi^{-1}(q)$ contains a further component
$F_{k+1}$ such that $F_1, \ldots, F_{k+1}$ are all components of $D_\ext^\vee$ with mother component $C^\vee_{t^\vee}$ contained in the fiber $\psi^{-1}(q)$. 
The components $F_1, \ldots, F_{k+1}$ form a chain in some feather $\fF^\vee$ of $D^{\vee\ge t^\vee}_\ext$.
\enum \elem

\bproof
All base points $p_{F_i}$ are equal since
by assumption the curves $F_i$ are contained in the same fiber over $q$. Let us denote this common base point by $p\in C_t$.

To deduce (a) we perform contractions in $D_\ext^{\ge t}$
until one of the components $F_0,\ldots, F_k$, say $F_0$, meets the
first time the component $C_t$.  On the contracted surface, say, $W'$,  the images of $F_0,\ldots, F_{k}$ necessarily form a chain
$[[-1,(-2)_k]]$. Any curve in $W$ between $F_0$ and $F_j$ has 
also mother component $C_t$ and the same base point  $p$ whence $F_0,\ldots, F_k$ form as well a chain in $D_\ext$.   Being connected this chain  has to be contained in some feather $\fF$.

(b) After renumbering we may suppose that any two consecutive curves in the chain $F_0,\ldots, F_k$ meet. The dual graph of $\fF_0=F_0+\ldots+ F_{k}$
in $W'$ is $[[-1,(-2)_k]]$, where $W'$ is as in (a). Since under the map $W\to W'$ only components of the fiber $\psi^{-1}(q)$ are contracted, the $\PP^1$-fibration $\psi:W\to\PP^1$ factors through a $\PP^1$-fibration $\psi':W'\to\PP^1$. The component $F_0$ meets $C_t$ in $W'$ and so  it disconnects the rest of the fiber $\psi'^{-1}(q)$ from $C_t$.
Hence this fiber cannot contain any component of the zigzag $D^{\ge t}$.

Being contractible in $W'$ to a smooth point $p \in C_t$
the chain $\cF_0$ cannot
exhaust the full fiber $\psi'^{-1}(q)$. 
After contracting $\fF_0$ in $W'$ the section $C_t$ 
still meets the resulting fiber in one point 
transversally. Hence there is an extra component 
$F_{k+1}$ of the fiber $\psi^{\prime -1}(q)$ with 
$F_{k+1}\cdot \fF_0=F_{k+1}\cdot F_k=1$.
Clearly, $F_{k+1}$ has
multiplicity $1$ in the fiber. Therefore the divisor $\psi'^{-1}(q)-F_k-F_{k+1}$ can be contracted on $W'$ to
a smooth point. Thus we arrive at a smooth surface $W''$ still fibered over $\PP^1$ with
two $(-1)$-curves $F_k$ and $F_{k+1}$ in the fiber over $q$.

The fiber component $F_{k+1}$ when considered as a curve on $W$,  cannot belong to $D^{\vee\ge t^\vee }$ since otherwise the feather $\fF$ containing $\fF_0$ would meet the boundary twice.
If on $W$ the fiber $\psi^{-1}(q)$ contains an extra component $C$ as in Lemma \ref{mem}(b), then $C$ must be contracted on $W''$.
Indeed, $C$ has to meet both $D^{\ge t}$ and $D^{\vee\ge t^\vee}$, whereas $F_{k+1}$, $F_k$ when considered in $W$ belong neither to $D^{\vee\ge t^\vee }$ nor to  
$D^{\ge t }$.

If $F_{k+1}$ were a component of
 $D_\ext^{\ge t}$ on $W$ then it would also be a
feather component with mother component $C_t$ and base point
$p$, contradicting the maximality of the collection $\{F_0,\ldots, F_k\}$.

According to Lemma \ref{mem}(b) $F_{k+1}$ is a component
of $D^{\vee\ge t^\vee}_\ext$ with
mother component $C^\vee_{t^\vee}$ and
base point $p^\vee\in C^\vee_{t^\vee}$ over $q$; this can be seen  from  the fiber structure on the surface $W''$.
Thus $F_{k+1}$ must be a component of a feather, say, $\fF^\vee$ of
$D_\ext^{\vee\ge t^\vee}$.

Assume that $k>0$. Then the feather $\fF$ containing $\fF_0$ can be written as $\fF=F+\fR$, where $F$ is the bridge curve of $\fF$ and $\fR$ contracts to a singular
point $x$ on $V$. The image of $F_{k+1}$ on $V$,
and then also that of $\fF^\vee$, contains $x$. It follows that
$\fF^\vee=F^\vee+\fR$ with $F^\vee$ being the bridge curve of $\fF^\vee$. In particular,
the chain $\fF_0-F_0+F_{k+1}=F_1+\ldots+F_{k+1}$ is contained
in $\fF^\vee$.

On the surface $W'$ we can contract all components of the fiber
$\psi^{\prime-1}(q)$
except for $F_0,\ldots,F_{k+1}$. The remaining fiber
$F_0+\ldots+F_k+F_{k+1}$
has then dual graph
$[[-1,(-2)_k,-1]]$. The subchain $F_1+\ldots+F_{k+1}$ of this fiber can be contracted on the resulting surface
to the point $p^\vee\in C^\vee_{t^\vee}$. Hence this
is a subchain of a maximal chain of feather components of
$D_\ext^{\vee\ge t^\vee}$ with the same mother component $C^\vee_{t^\vee}$
and the same base point $p^\vee$ as $F_{k+1}$.

If there were a further component $F^\vee$ of $D_\ext^{\vee\ge
t^\vee}$ with the same mother component
and the same base point as $F_{k+1}$,
then interchanging the roles of $\bV$ and
$\bV^\vee$ the above reasoning would give at least
$k+2$ feather components
 in $D_\ext^{\ge t}$ with mother component $C_t$ and base point
$q=p_{F_0}$, contradicting our assumption. Now the proof is
completed. 
\eproof

\bdefi\la{matching-pair} 
Following \cite[3.3.7]{FKZ4} a pair of feathers
$(\fF,\fF^\vee)$
as in Lemma \ref{matchlem} will be called a {\em matching pair}.
\edefi

The next lemma is shown in \cite[Lemma 3.3.10]{FKZ4} for a smooth
Gizatullin surface; the proof applies as well
to an arbitrary Gizatullin surface.

\blem \label{2.144} $C_t$ is a $*$-component if and only if
$C^\vee_{t^\vee}$ is. Furthermore, in the latter case the base points $q$
and $q^\vee$ as in Lemma \ref{mem}(a) coincide.
\elem

Now we are ready to deduce Proposition \ref{matchprop}.

\bsit\label{1.91} \bproof[Proof of Proposition \ref{matchprop}]
(a) is just Lemma \ref{2.144}. By Lemma \ref{matchlem}(b),  the
identification $C_t\cong C^\vee_{t^\vee}$ given by the
correspondence fibration $\psi$ provides a one-to-one
correspondence between the set of feathers of $D_\ext$ with mother
component $C_t$ and base point $p\in C_t$ and the set of feathers
of $D_\ext^{\vee}$ with mother component $C_{t^\vee}^\vee$ and  base point $p^\vee\in C_{t^\vee}^\vee$. Moreover,  the points
$$
\{\infty_{C_t}\}=C_t\cap C_{t-1}\and C^\vee_{t^\vee}\cap C_{t-1}=C^\vee_{t^\vee}\cap C^\vee_{t^\vee-1}=\{\infty_{C^\vee_{t^\vee}}\}\,
$$
are identified  under $\psi$ and as well the points $p$ and $p^\vee$ are. 
This leads  to the equality $0_{C_t}=0_{C^\vee_{t^\vee}}$ and  shows (b).
\eproof 
\esit

\subsection{Deformation invariance of the normalized extended graph}\label{definv}
The `only if' part of Theorem \ref{main1} is an immediate
consequence of the following proposition.

\bprop\label{order.70} If a family $\pi:\cV\to S$ of
affine ruled surfaces over a connected base $S$ admits a resolved
standard completion $(\tcV,\cD)$,
then the
normalized extended graph $N(\Gamma_\ext(s))$ of
$(\bar\cV(s),\bD(s))$ does not depend on $s\in S$. \eprop

\bproof The proof is based on the Factorization Lemma
\ref{factlem}. Given a  point $s_0\in S$
the fiber $(\tV,D)=(\tcV(s_0),\cD(s_0))$ of the pair $(\tcV,\cD)$ over $s_0$ is a resolved standard
completion of  $V$ with dual graph $\G$ independent of $s$. According to
Lemma \ref{factlem}, shrinking $S$ to a suitable \'etale
neighborhood of $s_0$ we
can decompose $\tcV\to S$ into a sequence
\be\la{fact1}
\tcV=\cX_m\to \cX_{m-1}\to \ldots \to \cX_1\to\cB\to S\,, 
\ee
where $\cB\to S$ is a smooth family of curves, 
$\cX_1\to\cB$ is a locally trivial $\PP^1$-fibration, 
and $\cX_i\to \cX_{i-1}$ is a blowup of a section 
$\gamma_i:S\to \cX_{i-1}$ with exceptional set 
$\cE_i\subseteq \cX_i$, $i=2,\ldots,m$. 
Restricting (\ref{fact1}) to the fiber
over $s_0$ we obtain a decomposition 
\be\la{fact2}
\tV=X_m\to X_{m-1}\to \ldots \to X_1\to \bB\,.
\ee 
Let now $F$ be a component of a feather $\fF$ of 
$D_\ext (s_0)$ created in the
blowup $X_i\to X_{i-1}$ with center on the mother component
$C\subseteq D$ of $F$. Then $\gamma_i(s_0)\in\cD$ is a smooth point of $\cD$ and so
$\gamma_i(S)\subseteq\cC$, where $\cC$ is a component of $\cD$ near $s_0$ with $\cC(s_0)=C$.
For $s\in S$ we let $F(s)$ be the proper transform of $\cE_i(s)$
in $\tcV(s)$ so that $F(s_0)=F$. As follows from the definition of a constant
dual graph of a family of completions (see Definition \ref{reldiv})
$F(s)$ cannot be a component of $\cD(s)$. 
Hence this is a feather component in $\tcV(s)$ 
with mother component $\cC(s)$. Since
every feather of $\tcV(s)$ appears in this way and the total number
of components of the extended graph $D_\ext(s)$ stays constant,
the number $\delta_C$
of feather components in $\tcV(s)$ with mother component $\cC(s)$ does not
depend on $s\in S$. Hence the normalized extended graph
$N(\Gamma_\ext(s))$ is indeed independent of $s\in S$. 
\eproof

The above proof provides in fact a little piece more of information.

\bcor\la{naturalcor}
The identification $N(\G_\ext(s_0))\cong N(\G_\ext(s))$ as
in the proof of Proposition \ref{order.70} is independent of
the factorization (\ref{fact1}). Thus locally in $S$ with respect to the \'etale topology there is a
{\em natural} identification $N(\G_\ext(s_0))\cong N(\G_\ext(s))$ for $s\in S$.
\ecor

\bproof
Given another factorization
\be\la{fact1a}
\tcV=\cX_m\to \cX'_{m-1}\to \ldots \to \cX'_1\to\cB\to S\,,
\ee
in the first step $\cX_m\to \cX_{m-1}$ and $\cX_m\to \cX'_{m-1}$
two families  of $(-1)$-curves $\cE$ and $\cE'$, respectively, are contracted. These relative $(-1)$-curves are either equal, or they are disjoint in all fibers by the uniqueness part in Lemma \ref{deflem}.  In the latter case it is possible to contract $\cE$ and $\cE'$ simultaneously.
In this way we obtain a new family $\tcV_1$ with again
two different factorizations. Using induction, these
two factorizations provide the same identification of extended divisors.
\eproof

The monodromy of the base points of the feather components can be
non-trivial  even in a family with constant dual graph.
To exclude this possibility we give the following definition.

\bdefi\la{naturaldef} 
Let $\cV\to S$ be a family of affine ruled surfaces
with resolved completion $(\tcV,\cD)$ and associated dual graph $\G$. We call
$(\tcV,\cD)$ a {\em family of constant type} $(\Delta, \G)$ if there is a family of
isomorphisms $\Delta\cong \G_\ext(s)$ such that for every point
$s_0\in S$ the isomorphism $N(\G_\ext(s))\cong \Delta\cong
N(\G_\ext(s_0))$ is the natural identification of Corollary
\ref{naturalcor} for $s$ near $s_0$. 
\edefi

\section{Versal families of affine ruled surfaces}\label{sec-versal}

\subsection{Complete deformation families of surfaces}

In this subsection we
construct a sufficiently big family of affine ruled surfaces admitting a
resolved completion with a given normalized extended graph $\Delta$.
This family is complete
in the sense that every other
such family can be obtained, at least locally,  from this one by
a base change. In particular, every individual affine ruled surface
admitting a resolved completion with normalized extended graph $\Delta$ appears
as a member of our family.

\bsit\la{3.1} 
We let $\G$ be a semi-standard tree with an extremal
$0$-vertex and a fixed embedding $\G\hookrightarrow\Delta$  into a
normalized extended weighted tree. All extremal 0-vertices say,
$v_{01},\ldots, v_{0a}$  of $\G$ are joined to the same vertex $v_1$,
which is uniquely determined. The subgraph $\Delta^1$ of $\Delta$
consisting of $v_1$ and all its neighbors in $\Delta$ will be
called the {\em one-skeleton} of $\Delta$. We consider it as a
weighted graph by assigning to $v_1$ the same weight as in $\Delta$,
while all other vertices get the weight zero.
Thus $\Delta^1$ consists of $v_1$, $v_{0i}$, $1\le i\le a$, and the
remaining neighbors say, $v_{21},\ldots,
v_{2b}$ of $v_1$. We assume that $\Delta$ is obtained from its
1-skeleton by a sequence of blow ups
\be\la{step-1} 
\Delta= \Delta^m\to \Delta^{m-1}\to\ldots \to
\Delta^1\,.
\ee
We suppose also that the following conditions are satisfied.

(a) If  in a blowup
$\Delta^i\to\Delta^{i-1}$ a vertex $v$ of $\G$ is created then in the subsequent
blowups we first create all the feathers of $v$, and only
after that the next vertices of $\G$.
In other words, if $v$ and $w$ are vertices of $\G$ created in the blowups
$\Delta^i\to\Delta^{i-1}$ and $\Delta^j\to \Delta^{j-1}$,
respectively, where $j>i$, then the vertices created in the
further blowups $\Delta^k\to \Delta^{k-1}$, $k>j$, are not
feathers of $v$.

(b) There is a genus $g$ assigned to $v_1$, and $v_1$ has weight $-2g$.
\esit

The construction of a versal family of surfaces associated to
$\Delta$ proceeds in several steps as follows.

\bsit\la{3.2}
 We let $\cM$ denote the moduli space of marked curves
$(C,p)$ of genus $g$ with $p\in C$ and with a level $l$ structure,
where $l\ge 3$. The latter means that we fix a symplectic basis of
the group $H^1(C;\ZZ/l\ZZ)$. It is known that $\cM$ is irreducible, quasi-projective 
\cite{DeMu} and  smooth \cite{Po}. According to Theorem 10.9 in
\cite[Lect. 10]{Po} there exists a universal family of curves
$\cZ\to \cM$, which constitutes a smooth projective morphism. 
\smallskip

\no{\bf Step 0:} By \cite[Th\'eor\`eme 3.2]{Gro0} the relative Picard functor
$\Pic_{\cZ/\cM}$ is representable over $\cM$. In view of the fact that there is  a section provided by the marking,  the representability of the Picard functor means that there exists  a
scheme $\cP$ locally of finite type over $\cM$ and a universal
line bundle $\cL$ over $\cZ\times_\cM\cP$, see \cite[Corollaire 2.4]{Gro0}. Letting
$\cP'\subseteq \cP$ be the connected component corresponding to the 
line bundles of degree $-2g$ and $\cB^0=\cZ\times_\cM\cP'$, there
are morphisms\footnote{By abuse of notation, we denote the
restriction $\cL|\cB^0$ still by $\cL$.} 
\be\la{step0}
\cX^0_{1}=\PP_{\cB^0}(\cO_{\cB^0} \oplus \cL)\to \cB^0\to
S^0=\cP'\,, 
\ee 
where the first one is a locally trivial
$\PP^1$-fibration and $\cB^0\to S^0$ is a complete family of
marked curves of genus $g$ with a level $l$ structure over $S^0$.
The composition $\cX^0_{1}\to S^0$ is a complete family of projective ruled
surfaces over $S^0$. The marking gives rise to a section
$\sigma:S^0\to\cB^0$. Taking the fiber product yields a family of
rational curves $\cC^0_{01}= \cX^0_{1}\times_\sigma S^0$  in
$\cX^0_{1}$ over $S^0$. The projection $\cO_{\cB^0} \oplus \cL\to
\cO_{\cB^0}$ defines a section $\cC_1^0\cong \cB^0$, where
$$
\cC_1^0=\PP_{\cB^0}(\cO_{\cB^0}) \hto
\cX^0_{1}=\PP_{\cB^0}(\cO_{\cB^0} \oplus \cL).
$$
Thus the fiber $X_1= \cX^0_{1}(t)$ over a point $t\in S^0$ contains the curves
$$
C_{01}= \cC_{01}^0(t)\quad \mbox{and}\quad C_1= \cC_1^0(t)\,.
$$
Clearly $X_1$ is a $\PP^1$-bundle over the genus $g$ curve
$B=\cB^0(t)$, the curve $C_{01}$ is a full fiber of $X_1\to B$ and
$C_1\cong B$ is a section. By construction $\cO_{C_1}(C_1)\cong
\cL|B$ so $C_1$ has  in
$X_1$ self-intersection $C_1^2=\deg (\cL|B)=-2g$.
\smallskip

\no{\bf Step 1:} Let now
$$
(\cB^0)^{a+b-1}=\cB^0\times_{S^0}\times\ldots\times_{S^0}\cB^0
$$
be the $(a+b-1)$-fold fiber product and $S^1\subseteq
(\cB^0)^{a+b-1}$ be the subset in the fiber over $s\in S^0$ consisting of all points
$$
(p_{02},\ldots ,p_{0a}, p_{21},\ldots , p_{2b})\,
$$
with pairwise distinct coordinates different from 
$p_{01}=\sigma(s)$. By base change $S^1\to S^0$ we 
obtain from (\ref{step0}) morphisms
\be\la{step1} 
\cX^1 _{1}\to\cB^1\to S^1\,. 
\ee 
On $\cX^1_{1}$ we have again the families of curves
$\cC_{01}^1=\cC_{01}^0\times_{S^0}S^1$ (corresponding to
the section $\sigma$) and
$\cC^1_1=\cC_1^0\times_{S^0}S^1$. With $\pi_i:S^1\to\cB^0$ being the
$i$th projection the morphisms
$$
\pi_i\times\id: S^1\to \cB^1\subseteq \cB^0\times_{S^0} S^1, \quad
i=1,\ldots, a+b-1,
$$
yield sections $\sigma_{02},\ldots, \sigma_{0a}$ and
$\sigma_{21},\ldots, \sigma_{2b}$. The preimages of
$\sigma_{ij}(S^1)\subseteq\cB^1$  under $\cX^1\to\cB^1$ give rise to
families of curves $\cC^1_{02},\ldots, \cC^1_{0a}$ and
$\cC^1_{21},\ldots, \cC^1_{2b}$ in $\cX^1$.
\smallskip

\no{\bf Further steps:} We construct in \ref{indstep}-\ref{nela} below a sequence of
morphisms 
\be\label{added} 
S^m\to S^{m-1}\to \ldots \to S^1\,,
\ee
which corresponds to the sequence (\ref{step-1}). 
Furthermore, for each $i=1,\ldots,m$ we define \bnum[(a)]
\item a sequence of morphisms
\be\la{stepi} 
\cX^i_i\to \cX^i_{i-1}\to\ldots \to \cX^i _{1}\to
\cB^i\to S^i, 
\ee
where $\cB^i=\cB^{i-1} \times_{S^{i-1}}S^i$, and
\item families of curves over $S^i$
$$
\cE^i_j\subseteq \cX^i, \,j\le i,\quad \cC^i_{0\alpha}\subseteq
\cX^i\,, \,\, 1\le \alpha\le a, \quad \cC^i_{1}\subseteq
\cX^i\,\mbox{ and } \, \cC^i_{2\beta}\subseteq \cX^i\,,\,\,
1\le\beta\le b,
$$
corresponding to the vertices of the 1-skeleton of $\Delta^i$ in  (\ref{step-1}), \enum
with the following
properties. 
\bnum[(i)]
\item Except for the morphism $\cX^i_i\to \cX^i_{i-1}$, (\ref{stepi})
is obtained from the
corresponding sequence of morphisms in level $i-1$ by base change
$S^i\to S^{i-1}$.
\item $\cE^i_i$ is a family of $(-1)$-curves obtained by blowing up a
section of $\cX^{i-1}_{i-1}\to S^{i}$. Moreover we have
$\cE^i_j=\cE^j_j\times_{S^j}S^i$ for  $j< i$ and similarly
$$
\cC^i_{0\alpha}=\cC^1_{0\alpha}\times _{S^1} S^i\,,\quad
\cC^i_{1}=\cC^1_{1}\times _{S^1} S^i \,\mbox{ and } \,
\cC^i_{2\beta}=\cC^1_{2\beta}\times _{S^1} S^i\,.
$$
\item The divisor $\cD^i_\ext$ on $\cX^i$ formed of the families of curves
$\cE^i_j$, $\cC^i_{0\alpha}$, $\cC^i_{1}$ and $\cC^i_{2\beta}$, restricts
in each fiber over $S^i$ to an SNC divisor with normalized dual graph
$\Delta^i$ as in (\ref{step-1}).
\item The divisor, say, $\cD^i$ formed of
families of curves as in (ii) that correspond to
vertices in $\G$, represents a family of SNC divisors with
dual graph $\G\cap \Delta^i$. \enum \esit

\bsit\la{indstep} Suppose that the sequences
(\ref{added}) and (\ref{stepi}) are already constructed up to a step $i<m$. Then at the
step $i+1$ we construct these data as follows.

(A) Suppose that in the blowup $\Delta^{i+1}\to\Delta^i$ in
(\ref{step-1}) a feather $v$ is created by an outer blowup of,
say, $v'\in \G\cap \Delta^i$. Let $\cC'$ be the corresponding
family of curves on $\cX^i$ so that $\cC'=\cE_j^i$ for some $j\le
i$, or $\cC'$ is one of the families of curves $\cC^i_{2,\beta}$
for some $\beta$. Letting
$$
S^{i+1}=\cC'\backslash \cD'\quad\mbox{with}\quad \cD'=\cD^i-\cC',
$$
we define $\cX_{i+1}^{i+1}$ to be the blowup of
$\cX_i^{i+1}=\cX_i^i\times_{S^i} S^{i+1}$ along $S^{i+1}$, which we
consider as a subscheme of $\cX_i^{i+1}$ via the diagonal embedding
into $\cX_i^i\times_{S^i} S^{i+1}$. Now we let $\cE^{i+1}_{i+1}$ be
the resulting family of $(-1)$-curves on $\cX_{i+1}^{i+1}$, while
the remaining curves in (b) are defined as in (ii).

(B) Suppose next that in the blowup $\Delta^{i+1}\to\Delta^i$
in (\ref{step-1})
a vertex $v$ of $\G$ is created. Then we proceed as follows
according to whether $v$ has one or two mother components.

(B1) If $v$ has one mother component, then it is created by an
outer blowup of, say, $v'\in \G\cap \Delta^i$. Let $\cC'$ be the
corresponding component of $\cD^i$. We let $S^{i+1}=\cC'\backslash
\cD'$, where $\cD'=\cD^i-\cC'$, and define now $\cX_{i+1}^{i+1}$
and $\cE_{i+1}^{i+1}$ similarly as in (A).

(B2) If $v$ has two mother components, then it is created by an
inner blowup of an edge connecting two vertices $v',v''\in \G\cap
\Delta^i$. These vertices correspond to families of curves $\cC'$ and
$\cC''$ of $\cD^i$, respectively. In this case we let
$S^{i+1}=S^i$, and we let $\cX^{i+1}_{i+1}$ be the blowup of
$\cX^i$ along the section $\cC'\cap\cC''$. 
\esit

\bsit\label{nela} 
On the last step $m$ we arrive at a family of surfaces
$\cX_j=\cX_j^m$ over $S=S^m$ with morphisms 
\be\la{stepn} 
\cX_m\to\cX_{m-1}\to\ldots\to \cX_1\to \cB=\cB^m \to S\,. 
\ee 
Likewise, omitting the upper index $m$  we obtain 
families of curves $\cC_{0\alpha}$, $\cC_1$, $\cC_{2\beta}$ 
and $\cE_j$. These yield a relative SNC divisor 
$\cD$ over $S$ with dual graph $\G$, and in every fiber of 
$\cX_m\to S$ over $s\in S$ also an extended divisor
$\cD_\ext(s)$. The normalized dual graph of $\cD_\ext(s)$ is
$\Delta$ for each point $s\in S$, while the dual graph of
$\cD_\ext(s)$ might depend on $s$.\footnote{So $\cD_\ext$ is not, 
in general, a relative SNC divisor over $S$.} 
We arrive in this way at a family of
smooth quasi-projective surfaces $\cX_m\backslash \cD$ over $S$.
In a final step using Proposition \ref{fcontract}, we can contract all feather components in 
$\cX_m\backslash \cD$ that are complete curves in $\cX_m\backslash \cD$. Thus we obtain a flat family of normal affine surfaces $\cV\to S$ along with a resolved $(-2g)$-standard
completion $(\tcV,\cD)=(\cX_m,\cD)$. 
\esit

\bdefi\la{3.5}
Given a normalized extended tree $\Delta$ together with a subtree $\G$ as in \ref{3.1} we call
\bnum[$-$]
\item $S=S_{\Delta,\G}$ as in \ref{nela} the {\em presentation space of normal affine surfaces
of type $(\Delta,\G)$};
\item  $\cV=\cV_{\Delta,\G}\to S$ the {\em universal  family  over $S$ of normal
surfaces of type $(\Delta,\G)$};
\item $(\tcV,\cD)=(\tcV_{\Delta,\G},\cD_{\Delta,\G})$  the {\em universal
resolved $(-2g)$-standard completion of $\cV\to S$ of type $(\Delta,\G)$}. \enum
\edefi

\brem 1. The construction of the presentation space $S_{\Delta,\G}$
and the universal families $(\tcV_{\Delta,\G}, \cD_{\Delta,\G})$
depends {\em \`a priori} on the order of blowups in the sequence
(\ref{step-1}). However, the reader can easily check that
different orders satisfying (a) and (b) in \ref{3.1}
will result in canonically isomorphic
presentation spaces and families over them.

2. The graph $\G$ does not determine the normalized extended
graph $\Delta$ as in \ref{3.1}, even if we restrict
to smooth Gizatullin surfaces. For instance, in the zigzag
$[[0,0,-2, -1, -2]]$ the component $C_3$ is a $*$-component, while
in $[[0,0,-1,-2,-1]]$ all of them are $+$-components. However,
blowing up suitable feathers we can obtain from both zigzags the chain
$\G=[[0,0,-2,-2,-2]]$. In the case of the zigzag $[[0,0,-2,-1,-2]]$ the resulting surface is an affine pseudo-plane, see \cite{MM} or \cite{FZ}, while in the other case the result is a Danilov-Gizatullin surface.

3. In the same way, given a normalized tree,  we can construct an associated universal resolved $s$-standard completion provided that $s\le -2g$. This is easily seen from the proof. If however $s>-2g$, $\cB$ is a family of complete curve of genus $g$ and $ \cL$ is a line bundle on $\cB$ of degree $s$ in each fiber, then fiberwise there are non-trivial extensions  $0\to \cO_{\cB_s}\to \cG\to \cL_s\to 0$. Even worse,  such extensions cannot be organized into a reasonable moduli space, in general.
\erem

In the next result we show that the universal family from Definition \ref{3.5} is complete in the sense of deformation theory. 

\bprop \la{uniprop}
Let $\cV'\to S'$ be a family of affine ruled surfaces
admitting a resolved $(-2g)$-standard completion 
$(\tcV',\cD')$.Then locally in the \'etale topology of $S'$ there is a morphism $S'\to S=S_{\Delta,\G}$ such that
$(\tcV',\cD')$ can be obtained from the universal family $(\tcV,
\cD)$ via a base change $S'\to S=S_{\Delta,\G}$. 
\eprop

\bproof 
By Proposition \ref{relstd} there is a  family $\cB'\to S'$ of curves of genus $g$  and the morphism $\tilde\Pi':\tcV'\to \tB$ such that for every $s\in S$
the restriction of $\tilde\Pi'$ to the fiber over $s$ is the standard morphism.
By the Factorization Lemma \ref{factlem} locally there is a
factorization 
$$
\tcV'=\cX'_n\to \cX'_{n-1}\to\ldots\to \cX'_1\to \cB' \to S'\,,
$$
where $\cB'\to S'$ is a smooth family of curves of genus $g$. 
In each step we blow up a section $S\to \cX'_i$, with the order of
blowups corresponding to that in (\ref{step-1}).

Locally the family $q':\cB'\to S'$ admits a level
$l$ structure. The position of the family of curves $\cC'_{01}$ in $\tcV'$ gives rise to a section $\sigma':S'\to\cB'$ and thus to a marking. Hence it is obtained from the universal family
$\cZ\to\cM$ in \ref{3.2} by a base change $S'\to \cM$. 
The $\PP^1$-bundle $\phi':\cX'_1\to \cB'$  possesses a
section $\cC_1'$ and so $\cX'_1\cong \PP^1_{\cB'}(\cG')$ is the projective bundle associated to the 2-bundle
$\cG'=\phi_*(\cO_{\cX'_1}(\cC_1'))$. The exact sequence
$$
0\to \cO_{\cX'_1}\to  \cO_{\cX'_1}(\cC_1') \to
\cL':=\cO_{\cC'_1}(\cC_1') \to 0
$$
yields a sequence $0\to \cO_{\cB'}\to \cG'\to \cL'\to 0$. In every
fiber over $s'\in S'$ the curve  $C_1'=\cC_1(s')$ satisfies
$C_1^{\prime 2}=-2g$, see \ref{3.1}. Hence $\cL'$ is a family of line bundles of
degree $-2g$ that is  obtained from the universal family
$\cL$ from \ref{3.2} by a base change $S'\to \cP'=S^0$. Furthermore the sequence
$0\to \cO_{\cB'}\to \cG'\to \cL'\to 0$ can be regarded as an
element in 
\be\la{extgroup}
\Ext^1_{\cO_{\cB'}}(\cL',\cO_{\cB'})\cong 
H^1(\cB', \cL^{\prime\,\vee})\,. 
\ee 
Since $\deg\,\cL^{\prime\,\vee}=2g$ the sheaf
$R^1q'_*(\cL^{\prime\,\vee})$ vanishes. Hence locally in $S'$
the group (\ref{extgroup}) vanishes as well. In other words, there is a
splitting $\cG'\cong \cO_{\cB'} \oplus\cL'$. It follows that
$\cX_1'\cong \cX_1\times_{S^0}S'$ locally  in $S'$.

Assume now that, for all $j\le i$,
$\cX'_j$ is already obtained from $\cX^i_j$ by a base change
$f_i:S'\to S^i$. The morphism $\cX'_{i+1}\to \cX'_i$ is then
the blowup of a section, say, $\sigma: S'\to \cX'_i$ with
exceptional set $\cE'_{i+1}\subseteq \cX_{i+1}'$, see Lemma
\ref{factlem}. Assume that in the blowup $\Delta^{i+1}\to
\Delta^i$ a vertex $v$ is created. As in \ref{indstep} we
distinguish between the following cases (A) and (B).

In case (A) the vertex $v$ corresponds to a feather component of $\Delta$ and is created
by an outer blowup at a unique vertex $v'\in \G$. The section
$\sigma$ must be contained in the corresponding family of curves,
say, $\cC'\subseteq \cX'_i$, which is necessarily a component of the image, say,
$\cD_i'$ of $\cD'$ in $\cX'_i$. The section $\sigma$ cannot meet any component of $\cD'_i$ different from $\cC'$ since otherwise in a fiber over some point $s'\in S'$
a feather is created by an inner blowup of
$D_i'=\cD_i'(s')$. Thus this newly created feather
will divide the zigzag, which is impossible. 

Using the
construction of $S^{i+1}$ in (A), the section $\sigma$ induces a
morphism $f_{i+1}=(f_i,\sigma): S'\to S^{i+1}$ such that
$\cX'_{i+1}\cong\cX_{i+1}\times_{S^{i+1}}S'$, as desired. 
The argument in case (B) of \ref{indstep} 
(where the newly created vertex $v$ is a vertex of $\G$)
is similar and so we leave it to the reader.
\eproof

\bcor\label{tercor} 
Let $V'$ and $V''$ be two normal affine ruled surfaces
admitting resolved standard completions $(\tV', D')$ and
$(\tV'', D'')$, respectively, with the same dual graph $\G$, the
same normalized extended tree $\Delta$, and the same genus assigned
to the vertex $v_1$ of $\G$. Then $V'$ and $V''$ are deformation
equivalent. 
\ecor

\bproof
Performing elementary transformations in one of 
the extremal 0-vertices $C_{01}'$ and $C_{01}''$
we may suppose that the curves $C_1'$ and $C_1''$ have the same
self-intersection index $-2g$. 
Then the induced completions $(\bV', D')$ and $(\bV'', D'')$ 
of $V'$ and $V''$, respectively,
arise as fibers over certain points $s',s''\in S$ of the
corresponding complete family $(\cX_m, \cD)$ over $S$. 
Since the moduli space of curves $\cM$ is irreducible, 
also the base $S$ is. 
In particular, $V'$ and $V''$ are deformation equivalent, as
required. \eproof

\subsection{The map into the configuration space}

We let as before $\phi: \cV\to S$ be a family of affine ruled surfaces, which admits a
minimal resolved completion $(\tcV,\cD)$ of constant type $(\Delta,\G)$.

Given a point $s\in S$, we consider the completion
$(\bV,D)=(\bcV(s),\cD(s))$  of the affine surface
$V=\cV(s)$ and the configuration invariant
$$ 
Q (\bV, D) \in \fC= \prod_{C\subseteq D-C_0-C_{1}}\fC_{\delta_C}(C^*)
\,,
$$
see Definition \ref{config-sp}. The configuration space $\fC$ does not depend on
the choice of $s\in S$ and only depends on the tree
$\Delta=N(\G_\ext)$ and the subtree $\Gamma=\G_D$. Note that $\delta_C$ is just the number of feathers of $\Delta$ at the vertex $C$.

\bdefi\la{3.8} $\fC=\fC(\Delta,\G)$ is called the {\em
configuration space associated to} $(\Delta, \G)$. \edefi

With these notation and assumptions we have the following result.

\bprop\la{3.9} The map $S\to\fC$ assigning to $s\in S$ the
configuration invariant of the fiber $\cV(s)$ is a regular morphism.
\eprop

\bproof 
This follows immediately from the Factorization Lemma
\ref{factlem}. Indeed, let
$$
\bcV=\cX_n\to \cX_{n-1}\to \ldots \cX_1\to \bar\cB\to S
$$
be the factorization as in Lemma \ref{factlem}. In each step
$\cX_{i+1}\to\cX_i$  the center of the blowup is a section $S\hto
\cX_i$ of $\cX_i\to S$. Thus restricting to the fiber over $s$ the
center of blowup varies regularly with $s$. This readily
implies the result. 
\eproof

\section{Moduli spaces of special Gizatullin surfaces}\la{modspecgiz}

\subsection{Special Gizatullin surfaces.}
In this section we discuss the existence of a coarse moduli space of affine ruled surfaces.
This does not follow immediately from the construction in section \ref{sec-versal}.
For instance, the moduli space of the Danilov-Gizatullin surfaces $V_n$ consists of one point, while  even in simple cases as in Example \ref{jump},
where $n=3$, our construction leads to a multi-dimensional versal family.  However, we show that in certain cases
the moduli space does exist and can be cooked out using the deformation family.

Fixing a tree $\G$ with an extremal $0$-vertex
as in \ref{3.1} and a normalized extended tree $\Delta$,
we consider families $\cV\to S$ of normal affine
surfaces over an algebraic $\kk$-scheme $S$ with a resolved completion $(\tcV,\cD)$ of constant type
$(\Delta, \G)$.
Two such families $\cX\to
S$ and $\cX'\to S$ are isomorphic if there is an
$S$-isomorphism $\cX\stackrel{\cong}{\longrightarrow} \cX'$.

\bdefi\la{3.10} Given an algebraic variety $S$ we let $\ff(S)$ be the
set of isomorphism classes of families $\cX\to S$  which admit a resolved completion of constant type $(\Delta,\G)$. This yields a functor
$$
\ff=\ff_{\Delta,\G}: \An\lto \sets
$$
from the category of algebraic $\kk$-schemes into sets. For a morphism
$S'\to S$ the corresponding map $\ff(S)\to \ff(S')$ is given by
the fiber product.

Restricting to families of smooth surfaces $\cX\to S$ we obtain
a functor $\ff_s: \An\lto \sets$.
\edefi

In general $\ff$ is not a sheaf. This means that two families
$\cX\to S$ and $\cX'\to S$ that are locally in $S$ isomorphic, do
not need to be  $S$-isomorphic globally. On the other hand, a
representable functor is always a sheaf. Thus in order to study
representability it is necessary to consider the sheaf $\tilde\ff$ associated to $\ff$. We present below concrete classes of surfaces
for which $\tilde\ff$ has a fine moduli space, which we
denote by $\fM(\Delta,\G)$. Clearly then $\fM(\Delta,\G)$ will be
as well a coarse moduli space for $\ff$.

\bexa\la{DGmod} 
Consider a Danilov--Gizatullin surface
$V_n=\Sigma_d\backslash C$ of type $n$, where $C\subseteq
\Sigma_d$ is an ample section with $C^2=n$ in the Hirzebruch
surface $\Sigma_d$. According to \cite{DaGi} (see also \cite{CNR,
FKZ5}) $V_n$ only depends on $n$ and neither on $d$ nor on the
choice of the section $C$. The normalized extended graph
$\Delta_n$ of $V_n$ is

\be\label{DGext0} \Delta_n:\qquad
\cou{0}{C_0}\llin\cou{0}{C_1}\llin \cou{\!\!\!\!\!\!\!\!\!\!\!\!
-2}{C_2}\nlin\cshiftup{F_1}{-1}
\llin\cou{-2}{C_{3}}\llin\ldots\llin\cou{\!\!\!\!\!\!\!\!\!\!\!\!
-2}{C_{n}} \nlin\cshiftup{F_0}{-1}\qquad. \ee 
By the Isomorphism
Theorem of Danilov and Gizatullin cited in the Introduction, the coarse moduli
space for such surfaces consists of a single reduced point.
\eexa

\bexa\label{spgis} 
According to \cite[Definition 1.0.4]{FKZ4} a {\em special Gizatullin
surface} with invariants $(n,r,t)$ is a smooth Gizatullin  surface with
normalized extended graph $\Delta=\Delta(n,r,t)$:

\vskip 0.1in
 \be\label{DGext1} \qquad
\cou{0}{C_0}\llin\cou{0}{C_1}\llin \cou{\!\!\!\!\!\!\!\!\!\!\!\!
-2}{C_2}\nlin\cshiftup{F_1}{-1}
\llin\cou{-2}{C_{3}}\llin\ldots\llin\cou{-2}{C_{t-1}}\llin
\cou{-2-r}{C_t}\nlin
\xbshiftup{}{\quad\;\;\{F_{ti}\}_{i=1}^r}\llin \cou{-2}{C_{t+1}}
\lin\ldots \lin \cou{\!\!\!\!\!\!\!\!\!\!\!\! -2}{C_{n}}
\nlin\cshiftup{F_0}{-1}\qquad, 
\ee
\normalsize
where $n\ge 3$ and
$\{F_{ti}\}_{i=1}^r$ is a family of $r$ feathers joined to $C_t$
and consisting each one of a single $(-1)$-curve. In the case
where $t=2$ or $t=n$ the number of  $(-1)$-feathers attached to
$C_t$ is $r+1$. Thus there are $\delta_t$ feathers with mother component $C_t$,
where $\delta_t=r+1$ for $t\in \{ 2, n\}$ and $\delta_t=r$ otherwise. The
associated configuration space is $\fC=\fC_{\delta_t}(C_t^*)$.

Assigning to $S$ the isomorphism classes of completable families over $S$
of special Gizatullin surfaces with invariants $(n,r,t)$ we obtain
as before a moduli functor $\bf F$. With this notation we can
reformulate Corollary 6.1.4 in \cite{FKZ4} as follows.
\eexa

\bthm\la{moduli} $\fC:=\fC_{\delta_t}(C_t^*)$ is a coarse moduli space
for $\bf F$. \ethm

\bproof By Proposition \ref{3.9} there is a functorial morphism
$$
\alpha_S: \ff(S)\to \Hom(S, \fC)\,.
$$
As follows from
Corollary 6.1.4 in \cite{FKZ4}, the elements of ${\bf F}(0)$ are
in one-to-one correspondence with the elements of $\fC$, where $0$
denotes the reduced point. It remains to show that for every other
space $\fM$ together with a functorial morphism 
$\beta: \ff(S)\to
\Hom (S,\fM)$ there is a unique morphism $\varphi:\fC\to\fM$ such
that the  diagram
\be\la{dia moduli}
\bdi
&&\ff(S)\\
&  \ldTo<\alpha &&\rdTo>{\beta}\\
\Hom(S,\fC)&&\rTo^{\varphi_*} &&\Hom(S,\fM)
\edi
\ee
commutes. Let $S(\G)$ be the space of
presentations as in \ref{3.5} and $\cV(\G)\to S(\G)$ be the
universal family.
Using $\beta$ this family induces a morphism
$$
\psi =\beta([\cV(\G)\to S(\G)]): S(\G)\to \fM.
$$
By Corollary 6.1.4 in \cite{FKZ4} this morphism is constant on the
fibers of $S(\G)\to \fC$. For any two elements in the fiber
define the same element in $\ff(0)$. Using the fact that $\fC$ is normal and $S(\G)\to \fC$ is surjective $\psi$ induces a
morphism $\varphi: \fC\to \fM$ making the diagram \eqref{dia moduli}
commutative, as required.
\eproof

\subsection{Appendix: An example} \la{sec-appendix}

Let us recall from Corollary \ref{semi-ample} that for a completable family of affine ruled surfaces $\cV\to S$ the total space $\cV$ is affine if so is $S$.

We emphasize that this does not remain true if we allow degenerations of $\cD$ in our completable families. Namely, allowing such degenerations we  show below the existence of a smooth family of affine surfaces $\cV\to S$ over $S=\AA^1$
with a completion
$(\bar\cV, \cD)$  of $\cV$ such that the total space $\cV$ is not affine and even not quasi-affine. We write in the sequel $\AA^n_{x_1,\ldots,x_n}=\Spec\kk[x_1,\ldots,x_n]$.

\bexa\la{exa2.2}
Similarly as in Example \ref{jump} we consider
the quadric $Q=\PP^1\times\PP^1$ with the zigzag $C_0+C_1+C_2$, where
$$
C_0=\{\infty\}\times\PP^1,\qquad C_1=\PP^1 \times\{\infty\},\quad\text{and}\quad C_2=\{0\}\times\PP^1\,.
$$
We let $V_1\to Q$ be the blowup in $(0,0)$ to create a feather $F$,  and $V_2\to V_1$ the blowup of $F\cap C_2$ to create a boundary curve $C_3$ on $V_2$. After choosing an isomorphism $j:S\cong C_2\backslash C_1$ we can blow up the image of the diagonal embedding $(1,j):S\hto S\times C_2$ on $S\times V_2$. Thus we obtain a surface $\bcV$ with a relative $(-1)$-curve $\cF'$ on it and a divisor $\cD=\sum_{i=0}^3\cC_i$, where $\cC_i$ is the proper transform of $S\times C_i$ in $\bcV$. We let $\cV$ be the complement $\bcV\backslash \cD$. The extended graphs of $(\bcV(s), \cD(s))$ over $s\ne 0$ and $s=0$ are 

\vskip2truemm
$$
\qquad \cou{0}{C_0}\lin\cou{0}{C_1}\lin
\cou{\qquad-2}{C_2}\nlin\cshiftup{F'}{-1}
\vlin{16}\cou{\qquad-2 }{C_3}\nlin\cshiftup{F}{-2}
\qquad\and \qquad \cou{0}{C_0}\lin\cou{0}{C_1}\lin
\cou{-2}{C_2}\lin\cou{-1}{F'}\lin\cou{\qquad-3 }{C_3}\nlin\cshiftup{F}{-2}\qquad ,
$$
respectively, where $F'$ stands for the fiber of $\cF'$ over $s$.
\eexa

\bprop\la{prop2.3}\bnum[(a)]
\item
Letting  $z=y/x$ we obtain
$$
H^0(\cV, \cO_\cV)\cong \kk[s,  zx, z^2x-sz,x ]\subseteq \kk[s,x,z]\,,
$$
while $H^0(\cV(0), \cO_{\cV(0)})=\kk[x, z]$.
\item
The variety $\cV$ is not quasi-affine.\enum
\eprop

\bproof
We only give a rough sketch of the argument; details will be given elsewhere. 

(a) implies (b), since the global functions on $\cV$ all vanish on the line $x=0$ of $\cV(0)$.

To deduce (a) let us simplify the setup by  considering  $V_i':=V_i\backslash (C_0\cup C_1)$, $i=1,2$, so that the restricted maps $V_2'\to V_1'\to \AA^2$ are blowups. In coordinates $V_1'\to \AA^2$
can be described by
$$
(x_1,\ y_1)=(x/y,\ y)\,\,,\mbox{where }
C_2 \mbox{ corresponds to the $y_1$-axis and
$F$ to the $x_1$-axis.}
$$
Clearly $V_1'\backslash C_2\cong \AA^2_{x,z}$ with $z=y/x$ since the other chart of the blowup $V_1'$ is given by $(x,z=y/x)$.
Blowing up $C_2\cap F$ with coordinates $(0,0)$ we obtain $V'_2$ with coordinates
$$
(x_2,\ y_2)=(x_1/y_1,\ y_1)=(x/y^2,\ y)\,,
$$
where
$$
C_2 \mbox{ corresponds to the $y_2$-axis and
the exceptional set $C_3$ to the $x_2$-axis.}
$$
The affine surfaces $V_2=V_2'\backslash (C_2\cup C_3)$ and $V_1=\AA^2_{x,z}$ are clearly isomorphic.
Next we have to blow up the product $S\times V_2'$ along the curve $\{(s,0,s)\,|\,s\in S\}$, where the last two coordinates are the $(x_2,y_2)$-coordinates as above. The ideal $I$ of this curve is given by 
$I=(x_2\,,y_2-s)$.
Blowing it up yields a 3-fold $\cV'$ with a coordinate chart $U_3\cong \AA^3$ and coordinates
\be\la{descr blowup}
(s,x_3,\ y_3)=(s, x_2,\ (y_2-s)/x_2)=(s,x/y^2, (y^3-sy^2)/x)\,,
\ee
where the new exceptional set $\cF'$ corresponds to $\{x_3=0\}$. By construction $\cV=\cV'\backslash(\cC_2\cup\cC_3)$. Moreover $\cC_3$ is given on $U_3$ by $y_2=0$, or in $(x_3,y_3)$ coordinates as $x_3y_3+s=0$.
The threefold $\cV$ is  equal to the union
of the two coordinate charts
$$
S\times \AA^2_{xz}\qquad\mbox{and}\qquad D(x_3y_3+s)\subseteq S\times  \AA^2_{x_3,y_3}\,.
$$
Accordingly
$$
A:=H^0(\cV,\cO_\cV)=\kk[s,x,z]\cap \kk[s,x_3,y_3, (x_3y_3+s)^{-1}]\,.
$$
Using \eqref{descr blowup} it is easy to see that  $A$ contains the functions
$$
s, \quad zx=x_3y_3+s, \quad x=\frac{(zx)^2}{z^2x},\quad z^2x-sz=\frac{z^3x^2-sz^2x}{zx} \,.
$$
It can be shown that $A$ is actually generated by them over $\kk$. Thus the result follows.
\eproof

\end{document}